\documentclass[10pt]{article}

\usepackage[english]{babel}
\usepackage{amsmath}
\usepackage{amsfonts}
\usepackage{amssymb}
\usepackage[all]{xy}
\usepackage{theorem}

\newcommand{\NN}{\mathbb{N}}
\newcommand{\CC}{\mathbb{C}}
\newcommand{\rra}{\rightrightarrows}
\newcommand{\ZZ}{\mathbb{Z}}
\newcommand{\HH}{\mathcal{H}}

\newcommand{\RR}{\mathbb{R}}
\newcommand{\RC}{\mathcal{R}}
\newcommand{\BB}{\mathcal{B}}

\newcommand{\la}{\left\langle}
\newcommand{\ra}{\right\rangle}

\newcommand{\fiber}[2]{\mathbin{{}_{#1}\times_{#2}}}
\newcommand{\inner}{\mathbin{\raise0.1ex\hbox{$\lrcorner$}}}

 \DeclareMathOperator{\supp}{supp}

 \DeclareMathOperator{\Hom}{Hom}

\DeclareMathOperator{\Res}{Res}

\DeclareMathOperator{\im}{im}

\DeclareMathOperator{\Rep}{\textbf{Rep}}

\DeclareMathOperator{\RRep}{Rep}
\DeclareMathOperator{\IdRep}{IdRep}
\DeclareMathOperator{\IrRep}{IrRep}

\DeclareMathOperator{\spann}{span}
\DeclareMathOperator{\Bis}{Bis}
\DeclareMathOperator{\PW}{PW}

\theoremstyle{plain}
\newtheorem{theorem}{Theorem}[section]
\newtheorem{lemma}[theorem]{Lemma}
\newtheorem{proposition}[theorem]{Proposition}

\newtheorem{corollary}[theorem]{Corollary}

{\theorembodyfont{\rmfamily}
\newtheorem{definition}[theorem]{Definition}
\newtheorem{remark}[theorem]{Remark}
\newtheorem{example}[theorem]{Example}
}

\newenvironment{proof}{\removelastskip\par\medskip
\noindent{\em Proof.}
\rm}{\penalty-20\null\hfill$\square$\par\medbreak}

\begin{document}

\title{Continuous representations of groupoids}

\author{Rogier Bos\thanks{Institute for Mathematics, Astrophysics and Particle Physics; Radboud University Nijmegen, The Netherlands, \texttt{Email:\ R.Bos@math.ru.nl},
this work is part of the research programme of the `Stichting voor Fundamenteel Onderzoek der Materie (FOM)', 
which is financially supported by the `Nederlandse Organisatie voor Wetenschappelijk Onderzoek (NWO)}}
\date{}
\maketitle
\begin{abstract}
We introduce unitary representations of continuous groupoids on continuous fields of Hilbert spaces.
We investigate some properties of these objects and discuss some of the standard constructions from representation theory in
this particular context. An important r\^{ole} is played by the regular representation.
We conclude by discussing some operator algebra associated to continuous representations of groupoids; in particular, we analyse the relationship of 
continuous representations of $G$ and continuous representations of the Banach $*$-category $\hat{L}^{1}(G)$.
\end{abstract}
\tableofcontents

\section*{Introduction}
Our purpose is to study some of the basic theory of continuous representations in the context of groupoids. 
Some work in this direction was initiated by Westman in \cite{wes}, \cite{wes2}. Representations of groupoids occur naturally in geometry, since the 
parallel transport associated to a flat connection on a vector bundle is a representation of the fundamental groupoid of the base space. Another 
place where they occur is as vector bundles over an orbifold, since these 
correspond to representations of the groupoid representing the orbifold. Also, for a group acting on a space, equivariant vector bundles over that 
space correspond to representations of the associated action groupoid.

We shall look at representations not only on continuous vector bundles, but on continuous fields of Hilbert spaces.
A reason why we not only consider representations on continuous vector bundles is the following.
One should note that the regular representation of a groupoid $G\rra M$ with Haar system is defined on a continuous field of $L^{2}$ functions on 
the target fibers. Even for very simple \'etale groupoids this is not a locally trivial field (consider the bundle of groups 
$(\ZZ/2\ZZ\times\RR)\backslash\{(-1,0)\}\to\RR)$.

We introduce the notion of a representation of a continuous groupoid on a continuous field of Hilbert spaces.
Continuous fields of Hilbert spaces were introduced and studied by Dixmier and Douady \cite{dido}. They play an important r\^ole  
in noncommutative geometry, as they occur as \mbox{(Hilbert $C^{*}$-)}modules of commutative $C^{*}$-algebras. 
Moreover, they are a rich source of noncommutative $C^{*}$-algebras, which are obtained as the algebra of adjointable endomorphisms of such modules.  

We develop an extension of harmonic analysis from
continuous groups to continuous groupoids.  It is investigated to which extent one can prove statements like Schur's Lemma and the Peter-Weyl 
theorem in the context of groupoids. Indeed, one can give an analogue of the decomposition of $L^{2}(G)$ for a compact group, under suitable 
conditions. We conclude by discussing some 
operator algebra associated to continuous representations of groupoids; in particular the relation to the continuous Banach $*$-category 
$\hat{L}^{1}(G)$.

Let us mention that representations of groupoids were also studied by J.\ Renault \cite{ren}. But one should note that the representations discussed there are measurable representations on measurable fields of Hilbert spaces. These behave quite different from continuous representations as studied in the present paper. 

As one will see in this paper, proofs of theorems in representation theory of groupoids heavily rely on the representation theory of groups. The differences mostly arise in dealing with the global topology of the groupoid and its orbit foliation.

Here follows an outline of the paper.

In the first section we resume some basic knowledge of continuous fields of Hibert spaces. This section contains no new material and its main purpose is to introduce and fix notation for continuous fields which are the main objects in the paper. 

The second section introduces representations of groupoids on continuous fields of Hilbert spaces. We discuss several notions of continuity
of representations and show how they relate. Then we treat two examples, namely the regular representation of a groupoid and representations of continuous families of groups. In the last part of the section we ``embed'' the theory of continuous groupoid representations in the theory of group representations. We discuss the topological group of global bisections of a groupoid and give a theorem that explains which representations of this group correspond to representations of the goupoid.

The third section treats harmonic analysis in the case of groupoids. We prove an analogue of Schur's Lemma and 2 versions of the Peter-Weyl Theorem. Then we consider Morita equivalence of groupoids and prove a theorem which states that Morita equivalent groupoids have equivalent representation categories.
The last part of this section discusses the representation rings of a groupoid.

The last section is a continuous analogue of Reneault's theorem that gives a bijection between measurable representations of $G$ and non-degenerate representations of the Banach algebra $L^{1}(G)$. We construct a bijection between continuous representations of $G$ and continuous non-degenerate representations of the Banach $*$-category $\hat{L}^{1}(G)$.

We shall denote a groupoid $G$ over $M$ by $G\rra M$. The source and target map are denoted by $s,t:G\to M$, the set of composable arrows
$G\fiber{t}{s}G$ by $G^{(2)}$, composition by $m:G^{(2)}\to G$, the unit map $u:M\to G$ and inversion by $i:G\to G$ or $g\mapsto g^{-1}$. 

The author would like to thank Gert Heckman, Peter Hochs, Klaas Landsman, Michael Mueger for helpful discussions on the topic and Maarten Solleveld
for some comments on a earlier version of this paper.

\section{Preliminaries: continuous fields.}
In this section we introduce continuous fields of Banach spaces and continuous fields of Hilbert spaces. We discuss the topology on the total space of 
such fields. A good understanding of this topology is crucial for many constructions in the rest of this paper. We then focus on uniformly 
finite-dimensional continous fields of Hilbert spaces and Lemma \ref{contrepstr} 
explains the structure of such a field. Finally, we discuss the relation of continuous fields of Banach/Hilbert spaces with Banach/Hilbert $C^{*}$
-modules. This is important since the regular representation of a groupoid is constructed from a specific Hilbert $C^{*}$-module. Most of the material 
in this section can be found \cite{dido}.  
\subsection{Continuous fields of Banach/Hilbert spaces.}
Let $M$ be a locally compact Hausdorff space. 
\begin{definition}\label{defcontf}
A \textbf{continuous field of Banach spaces over} $M$ is a family of Banach spaces
$\{\BB_{m}\}_{m\in M}$ and a space of sections
$\Delta\subset\prod_{m\in M}\BB_{m}$, such that
\begin{itemize}
\item[(i)] the set $\{\xi(m)\mid\xi\in\Delta\}$ equals $\BB_{m}$ for all $m\in M$.
\item[(ii)] For every $\xi\in\Delta$ the map $m\mapsto \|\xi(m)\|$ is in $C_{0}(M)$.
\item[(iii)] $\Delta$ is locally uniformly closed, i.e.\ if $\xi\in\prod_{m\in
M}\BB_{m}$ and for each $\varepsilon>0$ and
each $m\in M$, there is an $\eta\in\Delta$ such that
$\|\xi(n)-\eta(n)\|<\varepsilon$ on a neighborhood of $m$, then
$\xi\in\Delta$.
\end{itemize}
\end{definition}
There is a subclass of these continuous fields which has our special interest.
\begin{definition}
A \textbf{continuous field of Hilbert spaces over} $M$ is a family of Hilbert spaces
$\{\HH_{m}\}_{m\in M}$ over $M$ and a space of sections
$\Delta\subset\prod_{m\in M}\HH_{m}$ that form a continuous field of Banach spaces.
\end{definition}
\begin{remark}
The second condition in Definition \ref{defcontf} can then be replaced by the requirement that for any
$\xi,\eta\in\Delta$ the map
$m\mapsto\langle\xi(m),\eta(m)\rangle_{\mathcal{H}_{m}}$ is in $C_{0}(M)$.
The field is called \textbf{upper (lower) semi-continuous} if $m\mapsto \|\xi(m)\|$ is just upper (lower) continuous for every $\xi\in\Delta$. 
\end{remark}
\begin{lemma}\label{modul}
If $(\{\BB_{m}\}_{m\in M},\Delta)$ is a continuous field of Banach spaces, then $\Delta$ is a left $C_{0}(M)$-module. 
\end{lemma}
\begin{proof}
Let $f\in C_{0}(M)$ and $\xi\in\Delta$. Let $\varepsilon>0$ and $m\in M$ be given. Define
\[V_{m}:=\{n\in M\mid |f(m)-f(n)|<\frac{\varepsilon}{\|\xi(m)\|+1}\mbox{ and }|\|\xi(m)\|-\|\xi(n)\||<1\}\]
Then, for $n\in V$
\[\|f(n)\xi(n)-f(m)\xi(n)\|<\frac{\varepsilon}{\|\xi(m)\|+1}\|\xi(n)\|<\varepsilon.\]
Since $f(m)\xi\in\Delta$, we conclude by (iii) that $f\xi\in\Delta$.
\end{proof}
Actually $\Delta$ is a Banach $C^{*}$-module (see the next paragraph). 
\begin{lemma}
\label{toponfield}
If $(\{\BB_{m}\}_{m\in M},\Delta)$ is a continuous field of Banach spaces,
then there is a topology on the total space $\BB:=\coprod_{m\in M}\BB_{m}$ such that $\Delta$ 
equals the set of continuous sections $\Gamma_{0}(\BB)=C_{0}(M,\BB)$.
\end{lemma}
\begin{proof}
For each $\varepsilon>0$, $V\subset M$ open and $\xi\in\Delta$, we define
\[U(\varepsilon,\xi,V):=\{h\in\BB\mid\|h-\xi(p(h))\|<\varepsilon\mbox{ and }p(h)\in V\},\]
where $p:\BB\to M$ is the projection of the total space on the base.
One easily sees that these sets form a basis for a topology on $\BB$. Indeed, suppose that 
$U(\varepsilon_{1},\xi_{1},V_{1})$ and $U(\varepsilon_{2},\xi_{2},V_{2})$ are two of them and $h\in\BB$ lies in the intersection. 
By (1) there is a $\xi\in\Delta$ such that $\xi(m)=h$, where $m=p(h)$. Let $\varepsilon_{i}'=\varepsilon_{i}-\|h-\xi_{i}(m)\|$ for $i=1,2$.
Choose any $\varepsilon>0$ such that $\varepsilon<\varepsilon_{i}'$ for $i=1,2$.
Define
\[V:=\{m\in V_{1}\cap V_{2}\mid \|\xi(m)-\xi_{i}(m)\|<\varepsilon_{i}-\varepsilon\mbox{ for }i=1,2\}.\]
Then $U(\varepsilon,\xi,V)\subset U(\varepsilon_{1},\xi_{1},V_{1})\cap U(\varepsilon_{2},\xi_{2},V_{2})$.

Suppose $\xi\in\prod_{m\in M}\BB_{m}$ is a continuous section. Let $\varepsilon>0$ and $m\in M$ be given. Define $h:=\xi(m)$.
There is a $\xi'\in\Delta$ such that $\xi'(m)=h$. Let $V$ be any open neighborhood of $m$, then $W:=\xi^{-1}U(\varepsilon,\xi',V)$ is open and
on $W$ we have $\|\xi'-\xi\|<\varepsilon$. By (iii) we conclude that $\xi\in\Delta$.

Conversely, suppose $\xi\in\Delta$. Let $U(\varepsilon,\eta,V)$ be an open set in $\BB$, then
\[\begin{array}{ll}
\xi^{-1}U(\varepsilon,V,\eta)&=p(U(\varepsilon,\eta,V)\cap\xi(V))\\
&=\{m\in M\mid\|\xi(m)-\eta(m)\|<\varepsilon\}
\end{array}\]
Note that $\xi-\eta\in\Delta$, hence $m\mapsto\|\xi(m)-\eta(m)\|$ is continuous. We conclude that the above set is open, so that 
$\xi\in\Gamma_{0}(\BB)$.
\end{proof}
As a short notation we sometimes denote a continuous field of Banach spaces $(\{\BB_{m}\}_{m\in M},\Delta_{\BB})$ by $(\BB,\Delta)$. 
\begin{lemma}\label{normcont}
For any continuous field of Banach spaces $(\BB,\Delta)$ the map $\|.\|:\BB\to\RR_{\geq 0}$ is continuous.
\end{lemma}
\begin{proof}
Suppose $h\in\BB_{m}$ for certain $m\in M$. Given $\varepsilon>0$, take a $\xi\in\Delta$ such that $\xi(m)=h$ and
\[V:=\|\xi\|^{-1}(\|h\|-\varepsilon/2,\|h\|+\varepsilon/2).\]
This is an open set, since $\|\xi\|:M\to \RR_{\geq 0}$ is continuous.
So, $h'\in U(\varepsilon/2,\xi,V)$, with $h'\in\BB_{m'}$ implies 
\[|\|h'\|_{m'}-\|h\|_{m}|\leq\|h'-\xi(m')\|+|\|\xi(m')\|_{m'}-\|h\|_{m}|\leq\varepsilon,\]
which finishes the proof.
\end{proof}

\begin{definition}\label{morphcontf}
A \textbf{morphism $\Psi:(\BB^{1},\Delta^{1})\to (\BB^{2},\Delta^{2})$ of continuous fields of Banach spaces} is a family of linear maps 
$\{\Psi_{m}:\BB^{1}_{m}\to\BB^{2}_{m}\}_{m\in M}$ such that the induced map $\Psi:\BB_{1}\to\BB_{2}$ on the total spaces satisfies
\[\{\Psi\circ\xi\mid\xi\in\Delta^{1}\}\subset\Delta^{2}\]
and 
\[m\mapsto\|\Psi_{m}\|\]
is a locally bounded map.
\end{definition}\label{defmap}
Here $\|\Psi_{m}\|$ is the operator norm of $\Psi_{m}$,
\[\|\Psi_{m}\|:=\sup_{\|h\|_{\BB^{1}_{m}}=1}\|\Psi_{m}(h)\|_{\BB^{2}_{m}}.\]
The first condition has to be satisfied only on a dense subset of $\Delta^{1}$ (\cite{dido}, Proposition 5).
\begin{lemma}
The map $\Psi:\BB^{1}\to\BB^{2}$ is continuous iff $\Psi$ is a morphism of continuous fields of Banach spaces.
\end{lemma}
\begin{proof}
"$\Leftarrow$" Suppose $h\in U(\epsilon_{2},\xi_{2},V_{2})\subset\BB_{2}$ and $p(h)=m$. By (i), there is a $\xi_{1}\in\Delta_{1}$ such that
$\xi_{1}(m)=h$. Since $\Psi(\xi_{1})\in\Delta_{2}$, the set defined by
\[V_{1}:=\{n\in M\mid\|\Psi(\xi_{1})-\xi_{2}\|(n)<\varepsilon/2\}\cap V_{2}\]
is open. Let $f:M\to\RR$ be a locally bounded function such that $\|\Psi(\xi)\|<f\|\xi\|$ for all $\xi\in\Delta$.
Let $V_{1}'\subset V_{1}$ be a small enough neighborhood of $m$ such that $f$ has a supremum $K$ on $V_{1}'$, then
\[\Psi(U(\frac{\varepsilon_{2}}{2 K},\xi_{1},V_{1}')\subset U(\varepsilon_{2},\xi_{2},V_{2})\]
Indeed, for any $h'\in U(\frac{\varepsilon_{2}}{2 K},\xi_{1},V_{1}')$ with $p(h')=n$ we have
\[\begin{array}{ll}
\|\Psi(h')-\xi_{2}(n)\|&=\|\Psi(h')-\Psi(\xi_{1}(n))+\Psi(\xi_{1}(n))-\xi_{2}(n)\|\\
&\leq\|\Psi(h'-\xi_{1}(n))\|+\|\Psi(\xi_{1}(n))-\xi_{2}(n)\|\\
&=K\|h'-\xi_{1}(n)\|+\frac{\varepsilon_{2}}{2}\\
&=K \frac{\varepsilon_{2}}{2 K}+\frac{\varepsilon_{2}}{2}=\varepsilon_{2}.
\end{array}\]
"$\Rightarrow$"  $\Psi(\Delta_{1})\subset\Delta_{2}$ by Lemma \ref{toponfield}. Let $m\in M$ be any element. By continuity $\Psi^{-1}(U(1,0,M))$ is open,
so it contains an open neighborhood $U(\varepsilon,0,V)$, where $V$ is an open neighborhood of $m$. Hence, $\|\Psi\|$ is bounded on $V$.
\end{proof}

The map $\Psi:(\BB^{1},\Delta^{1})\to (\BB^{2},\Delta^{2})$ is an (isometric) isomorphism of continuous fields of Banach spaces if all the $\Psi_{m}$ are (isometric) isomorphisms and $\Psi(\Delta^{1})=\Delta^{2}$. In fact, 
one can replace the second condition by $\Psi(\Lambda)\subset\Delta^{2}$ for a dense subset $\Lambda\subset\Delta^{1}$ (\cite{dido}, Proposition 6).

Let $(\{\BB_{m}\}_{m\in M},\Delta)$ be a continuous field of Banach spaces over $M$ and $J:N\to M$ a continuous map.
Define the \textbf{pullback continuous field} $J^{*}(\{\BB_{m}\}_{m\in M},\Delta)$ as follows. The fiber at $n\in N$ is $\BB_{J(n)}$. The space of
section $J^{*}\Delta$ is the smallest Banach space generated by $f\cdot J^{*}\xi$ for all $\xi\in\Delta$ and $f\in C_{0}(N)$, such that one obtains a
continuous field of Banach spaces. It is the closure of $\{\xi\circ J\mid\xi\in\Delta\}$ as a Banach space (see next the paragraph).
The continuous field thus obtained is denoted by $(J^{*}\{\BB_{m}\}_{m\in M},J^{*}\Delta)$. We shall need this construction in Section \ref{morequ}.

\subsection{Uniformly finite-dimensional continuous fields of Hilbert spaces.}
The \textbf{dimension of a continuous field of Hilbert spaces} is the supremum of the dimensions of its fibers. Note that dimension is 
a lower semi-continuous functions $M\to\ZZ\subset\RR$. That is, $\dim:M\to\ZZ_{\geq 0}$ has a local minimum at every point.
A continuous field of Hilbert spaces is \textbf{uniformly finite-dimensional} if it has finite dimension. 
One should distinguish between uniformly finite-dimensional and \textbf{finite-dimensional} continuous fields, which means 
that each fiber is finite dimensional.
\begin{example}
Consider the field over $\RR$ with $\HH_{x}:=\CC^{n}$ if $x\in [-n,-n+1)\cup(n-1,n]$ for $n\in\ZZ_{\geq 0}$.
This field is finite-dimensional, but not uniformly finite-dimensional. 
\end{example}

The following lemma characterizes uniformly finite-dimensional continuous fields of Hilbert spaces. Suppose $d\in\mathbb{N}$. 
\begin{lemma}\label{contrepstr}
A $d$-dimensional continuous field of Hilbert spaces $(\HH,\Delta)$ over a locally compact Hausdorff space $M$ is a sum 
\[\HH\simeq\sum_{i\in I}E_{i}\]
of vector bundles $E_{i}\to U_{i}$, where $\{U_{i}\}_{i\in I}$ is a locally finite open covering of $M$.
\end{lemma}
\begin{proof}
For each $m\in M$ choose sections $\{\xi_{j}^{m}\}_{j=1}^{\dim(\HH_{m})}$, such that $\{\xi_{j}^{m}(m)\}_{j=1}^{\dim(\HH_{m})}$ forms a basis of $\HH_{m}$. Let $V_{m}$ be the
set on which their images stay linearly independent and non-zero. This set is open, since
\[m\mapsto\det(\xi_{1}^{m}\mid\ldots\mid,\xi_{\dim{\HH_{m}}}^{m})=\det((\xi_{k}^{m}\cdot\xi_{l}^{m})_{k l})\]
is a continuous map. Indeed, this last expression is a polynomial in $\xi_{k}^{m}\cdot\xi_{l}^{m}$ for $1\leq k,l\leq j$ which are in $C_{0}(M)$ by definition. Define a subfield by
\[E_{m}:=\spann\{\xi_{j}^{m}\mid m\in V_{m}, j=1,\ldots,\dim(\HH_{m})\}.\]
One easily sees that these are indeed vector bundles, since the $V_{m}$ are trivializing neighborhoods and
trivializing diffeomorphisms are
\[\left.H\right|_{V_{m}}\to U_{i}^{m}\times\CC^{i}, h\mapsto(p(h),\la\xi_{1}^{m},h\ra,\ldots,\la\xi_{i}^{m},h\ra),\]
which finishes the proof.
The collection $\{V_{m}\}_{m\in M}$ covers $M$, hence there is a locally finite subcover $\{U_{i}\}_{i\in I}$ with vector bundles $\{E_{i}\to U_{i}\}_{i\in I}$. These are the desired vector bundles.
\end{proof}
\begin{corollary}
A continuous field $(\HH,\Delta)$ over a compact space $M$ is uniformly finite-dimensional iff $\Delta$ is finitely generated over $C_{0}(M)$.
\end{corollary}
Uniformly finite-dimensional continuous fields of Hilbert spaces over compact spaces arise as the regular representation of families of
finite groups, cf. Section \ref{secregrep}.

\begin{definition}
A continuous field $(\HH,\Delta)$ is \textbf{locally trivial} if for every $m\in M$ there exist a neighborhood $U\ni m$, a Hilbert space $\HH'$ and an 
isomorphism of continuous fields $\HH|_{U}\to U\times \HH'$.
\end{definition} 
\begin{example}
Locally trivial finite-dimensional continuous fields of Hilbert spaces are vector bundles.
\end{example}
\begin{definition}\label{loctriv}
A \textbf{local trivialization} around $m$ of a continuous field $(\HH,\Delta)$ is a neighborhood $U\ni m$, a Hilbert space $\HH'$ of dimension 
$\dim(\HH|_{U})$ and an injective morphism of continuous fields $\HH|_{U}\to U\times \HH'$.
\end{definition}
Such trivializations are important in Section \ref{famgr}.

\subsection{Banach/Hilbert $C^{*}$-modules.}
Let $A$ be a $C^{*}$-algebra and $A^{+}$ the set of positive elements in $A$.
\begin{definition}
A \textbf{left Banach $A$-module} is a Banach space $\Delta$, which has a left $A$-module
structure $A\to \BB(\Delta)$ and a linear map
$\|\cdot\|:\Delta\to A^{+}$ such that for all
$\xi,\eta,\chi\in\Delta$ and $a\in A$:
\begin{itemize}
\item[(i)] the norm on $\Delta$ satisfies $\|\xi\|_{\Delta}=\sqrt{\|(\|\xi\|^{2})\|_{A}}$,
\item[(ii)] $\|\xi+\eta\|\leq\|\xi\|+\|\eta\|$,
\item[(iii)] $\|a \xi\|=|a|\|\xi\|$, where $|a|:=\sqrt{a^{*}a}$,
\item[(iv)] $\|\xi\|=0$ iff $\xi=0$.
\end{itemize}
\end{definition}
As in the case of continuous fields, one has the subclass of Hilbert $A$-modules.
\begin{definition}
A \textbf{left Hilbert $A$-module} is a Banach space $\Delta$, that has a left $A$-module
structure $A\to \BB(\Delta)$ and a sesquilinear pairing
$\la\cdot,\cdot\ra:\Delta\times\Delta\to A$ such that for all
$\xi,\eta,\chi\in\Delta$ and $a\in A$:
\begin{itemize}
\item[(i)] the norm on $\Delta$ satisfies $\|\xi\|_{\Delta}=\sqrt{\|\la\xi,\xi\ra\|_{A}}$,
\item[(ii)] $\la\xi,\eta+\chi\ra=\la\xi,\eta\ra+\la\xi,\chi\ra$,
\item[(iii)] $\la\xi,a\eta\ra=a\la\xi,\eta\ra$,
\item[(iv)] $\la\xi,\eta\ra=\la\eta,\xi\ra^{*}$,
\item[(v)] $\la\xi,\xi\ra>0$ iff $\xi\not=0$.
\end{itemize}
\end{definition}
A \textbf{morphism of Banach $C_{0}(M)$-modules} is an operator $\Psi:\Delta_{1}\to\Delta_{2}$, that intertwines the $C_{0}(M)$ action and 
is such that $\|\Psi\|$ is a locally bounded map $M\to\RR$.
\begin{theorem}
Then there is an equivalence of categories of continuous fields of Banach (resp. Hilbert) spaces and
left Banach (resp. Hilbert) $C_{0}(M)$-modules.
\end{theorem}
\begin{proof} (sketch, for a full proof see \cite{dido} \S4),
Suppose $(\HH,\Delta)$ is a continuous field of Hilbert spaces. Then $\Delta$ is a $C_{0}(M)$-module as proven in Lemma \ref{modul}. Its
completeness as a Banach space follows immediately from locally uniform completeness. This is one direction of the correspondence.

For the other direction, suppose $\Lambda$ is a Banach $C_{0}(M)$-module. Define
\[N_{m}:=\{h\in\Lambda\mid \|h\|(m)=0\}\]
and $\HH_{m}:=\Lambda/N_{m}$. Denote the projection by $\pi_{m}:\Lambda\to\Lambda/N_{m}$.
Define the space of sections by
\[\Delta:=\{\xi_{\lambda}:=(m\mapsto\pi_{m}(\lambda))\mid \lambda\in\Lambda\}.\]
We check that this is indeed a continuous field of Hilbert spaces.
\begin{itemize}
\item[i)] $\{\xi_{\lambda}(m)\mid\xi_{\lambda}\in\Delta\}=\Lambda/N_{m}$ trivially;
\item[ii)] $m\mapsto\|\xi_{\lambda}(m)\|=\|\lambda\|(m)$ is by definition continuous;
\item[iii)] suppose $\lambda\in\prod_{m\in M}\Lambda/N_{m}$ and suppose $\lambda$ is locally uniformly close to sections in $\Delta$. We want to
show that this implies $\lambda\in\Delta$. Since $\Lambda$ is complete as a Banach space it suffices to show globally uniformly close to a
section in $\Delta$. This one shows using a partition of unity argument. We omit the details.
\end{itemize}

If one begins with a Banach $C_{0}(M)$-module $\Lambda$, then produces a continuous field of Banach spaces, and from that again constructs a Banach
$C_{0}(M)$-module, one trivially recovers $\Lambda$.

On the other hand, from a continuous field $(\{\HH_{m}\}_{m\in M},\Delta)$ one obtains the Banach $C_{0}(M)$-module $\Delta$ and once again this gives
rise to a continuous field $(\{\Delta/N_{m}\}_{m\in M},\Delta)$. An isomorphism $\Delta/N_{m}\to \HH_{m}$ is given by $[\xi]\mapsto\xi(m)$.
\end{proof}
\begin{remark}
This correspondence shall be extended in Theorem \ref{bisect}.
\end{remark}
The well-known Serre-Swan theorem states that for compact $M$ there exists an equivalence of categories between finitely generated projective
Hilbert $C(M)$-modules and locally trivial finite-dimensional continuous fields of Hilbert spaces (i.e. finite rank vector bundles) over $M$. Indeed,
as mentioned on compact spaces $M$ finitely generated Hilbert modules $\Delta$ correspond to uniformly finite-dimensional continuous fields. Moreover, one can show that $\Delta$ being projective corresponds to the field being locally trivial.

\begin{example}\label{exfield}
Suppose $\pi:N\to M$ is a continuous surjection. A continuous family of measures on $\pi:N\to M$ is a family of measures $\{\nu_{m}\}_{m\in M}$ on $N$
such that the support of $\nu_{m}$ is in $\pi^{-1}(m)=:N_{m}$ and for every function $f\in C_{c}(N)$ the function
\[m\mapsto \int_{N_{m}}f(n)\,\nu_{m}(dn)\]
is continuous $M\to\CC$.

For any $p\in\RR_{\geq 1}$ consider the norm on $C_{c}(N)$ given by
\[\|f\|_{p}:=\sup_{m\in M}\|f|_{N_{m}}\|_{L^{p}(N_{m},\nu_{m})}.\]
Define $\Delta_{\pi}^{p}(N)$ to be the closure of $C_{c}(N)$ with respect to this norm.
One easily sees that this is a Banach $C_{0}(M)$-module with $C_{0}(M)$-valued norm given by
\[\|f\|(m):=\|f|_{N_{m}}\|_{L^{p}(N_{m},\nu_{m})}=\left(\int_{N_{m}}|f(n)|^{p}\nu_{m}(dn)\right)^{1/p}.\]

The continuous field associated to this Banach $C_{0}(M)$-module is denoted by $(\hat{L}^{p}_{\pi}(N),\Delta_{\pi}^{p}(N))$.
The fiber at $m\in M$ equals $L^{p}(N_{m},\nu_{m})$.

If $p=2$ one obtains a Hilbert $C_{0}(M)$-module and hence a continuous field of Hilbert spaces. The $C_{0}(M)$-valued inner product is given on
$C_{c}(N)$ by
\[\la f,f'\ra(m):=\la f|_{N_{m}}, f'|_{N_{m}}\ra_{L^{2}(N_{m},\nu_{m})}=\int_{N_{m}}\overline{f(n)}f'(n)\nu_{m}(dn).\]
\end{example}

\section{Continuous representations of groupoids}
\subsection{Representations of a groupoid on a continuous field of Hilbert spaces.}\label{defrep}
In this section we introduce continuous representations of groupoids on continuous fields of Hilbert spaces. As far as we know this notion as we 
define it does not appear anywhere in the literature. We should mention the work of Westman \cite{wes, wes2} though, who restricts himself to 
representations of locally trivial groupoids on vector bundles. Furthermore, there is a preprint by Amini \cite{amin}, which treats continuous 
representations on families of Hilbert spaces with some continuity condition, but without condition i) of Definition \ref{defcontf}.

As for representations of groups there are
several forms of continuity for such representations. We consider ``normal'' and weak, strong continuity and in Section \ref{opno} also continuity in 
the operator norm. All these forms of continuity can be compared, cf.\ Lemma \ref{strongweak}, Lemma \ref{normalstrong} and Lemma \ref{opnormal},
generalizing similar results for groups (cf.\ \cite{gaal}).
In Definition \ref{morphrep} we introduce morphisms of representations and we show in Proposition \ref{makeun} that any 
representation of a proper groupoid is isomorphic to a unitary representation, generalizing a similar result for compact groups.

Let $M$ be a locally compact space and $G\rra M$ a continuous groupoid. 
\begin{definition}\label{boundrep}
A \textbf{bounded representation} of $G$ is a continuous field of Hilbert spaces $(\{\HH_{m}\}_{m\in M},\Delta)$ over
$M$ and a family of invertible bounded operators 
\[\{\pi(g):\HH_{s(g)}\to\HH_{t(g)}\}_{g\in G}\]
satisfying
\begin{itemize} 
\item[i)] $\pi(1_{m})=id_{\HH_{m}}$ for all $m\in M$,
\item[ii)] $\pi(g g')=\pi(g) \pi(g')$ for all $(g,g')\in G^{(2)}$,
\item[iii)] $\pi(g^{-1})=\pi(g)^{-1}$ for all $g\in G$ and
\item[iv)] $g\mapsto\|\pi(g)\|$ is locally bounded.
\end{itemize}
\end{definition}
We denote such a representation by a triple $(\HH,\Delta,\pi)$.

\begin{definition}
A representation $(\HH,\Delta,\pi)$ is \textbf{strongly continuous} if the map 
\[g\mapsto \pi(g)\xi(s(g))\] 
is continuous $G\to\HH$ for all $\xi\in \Delta$.
A representation is \textbf{weakly continuous} if the map 
\[g\mapsto\la\pi(g)\xi(s(g)),\eta(t(g))\ra\] 
is continuous $G\to\CC$ for all $\xi,\eta\in \Delta$.
A representation $(\pi,\HH,\Delta)$ is \textbf{continuous} if
\[\Psi: (g,h)\mapsto\pi(g)h\]
is a continuous map $G\fiber{s}{p}\HH\to\HH$.
The representation is \textbf{unitary} if the operators \mbox{$\{\pi(g):\HH_{s(g)}\to\HH_{t(g)}\}_{g\in G}$} are unitary.
\end{definition}
For any $\xi,\eta\in\Delta^{\pi}$ we use the notation $\la\xi,\pi\eta\ra$ for the map $G\to\CC$ given by
\[g\mapsto\la\xi(t(g))),\pi(g)\eta(s(g))\ra.\]

Condition (iv) of Definition \ref{boundrep} is perhaps somewhat strange at first sight. The following Example \ref{exbound1}, Lemma \ref{lsc} and
Example \ref{exbound2} should clarify it. Moreover, recall that for morphism $\Psi$ of continuous fields the map $m\mapsto \|\Psi_{m}\|$ has to be locally bounded too, cf.\ Definition \ref{morphcontf}.
\begin{example}\label{exbound1}
A simple example shows that $g\mapsto\|\pi(g)\|$ is not always continuous. Consider the groupoid $\RR\rra\RR$, 
with a continuous representation on a field given by 
trivial representation on $\CC$ at each $x\in\RR$ except in $0$, where it is the zero representation. 
In this case the norm of $\pi$ drops from $1$ to $0$ at $0$.  
\end{example}

\begin{lemma}\label{lsc}
For any continuous representation $(\HH,\pi,\Delta)$ the map $g\mapsto\|\pi(g)\|$ is lower semi-continuous $G\to\RR$.
\end{lemma}
\begin{proof}
Using the above definition and Lemma \ref{normcont} we know that the map $(g,h)\mapsto\|\pi(g)h\|$
is continuous $G_{s}\times_{p}\HH\to\RR_{\geq 0}$. For any $g\in G$, let $\varepsilon>0$ be given.
Let $h'\in\HH_{s(g)}$ be such that
\[|\|\pi(g)h'\|-\|\pi(g)\| |<\varepsilon/2.\]
by continuity there exists an open neighborhood $U\subset G_{s}\times_{p}\HH$ of $(g,h')$ such that $(g'',h'')\in U$ implies
\[|\|\pi(g'')h''\|-\|\pi(g')h'\||<\varepsilon.\]
Take $V:=pr_{1}(U)\subset G$. Then $g''\in V$ implies, for any $h''\in pr_{2}(U)$,
\[\|\pi(g'')\|\geq\|\pi(g'')h''\|>\|\pi(g)\|-\varepsilon,\]
and we are done.
\end{proof}

The function $g\mapsto\|\pi(g)\|$ is locally bounded if, for example, $(\HH,\Delta)$ is finite-dimensional.
\begin{example}\label{exbound2}
A counterexample of a continuous representation of a proper groupoid where $g\mapsto\|\pi(g)\|$ is not locally bounded $G\to\RR$, even though 
the restriction to $G_{m}$ is bounded for each $m$, is as follows.

Consider the trivial bundle of groups $[0,1]\times \ZZ/2\ZZ\rra [0,1]$. Define a continuous field of Hilbert spaces over $[0,1]$ by
$\HH_{0}:=\CC^{2}=:\HH_{1}$ and $\HH_{x}:=\CC^{2\,n}$ if $x\in [\frac{1}{n+1},\frac{1}{n})$ for all $n\in\NN$. The topology on the field
is obtained from the inclusions $\CC^{2n}\hookrightarrow\CC^{2(n+1)}$ given by $\vec{v}\mapsto(0,\vec{v},0)$. Define, for every $n\in\NN$ and 
$x\in [\frac{1}{n+1},\frac{1}{n})$,
\[\pi(x,-1):=\mbox{diag'}(1/n,\ldots,1/2,1,1,2,\ldots n),\]
where diag' denotes the matrix filled with zeros except the diagonal from the upper right corner to the lower left corner, where the above sequence 
is filled in. Furthermore, $\pi(0,-1):=\mbox{diag'}(1,1)$. This representation is strongly continuous, 
but 
\[\|\pi(x,-1)\|=n\mbox{ if }x\in\left[\frac{1}{n+1},\frac{1}{n}\right).\]
Hence $g\mapsto\|\pi(g)\|$ is not locally bounded at $(0,-1)$.
\end{example}

\begin{lemma}
\label{strongweak}
If a representation $(\pi,\HH,\Delta)$ is strongly continuous, then it is weakly continuous. The converse implication holds if the 
representation is unitary.
\end{lemma}
\begin{proof}
Suppose $(\pi,\HH,\Delta)$ is strongly continuous.
Suppose $\xi,\eta\in\Delta^{\pi}$ and $g\in G$. Write $n=t(g)$. Let $\varepsilon>0$ be given. Let $\xi'\in\Delta^{\pi}$ be a section satisfying
$\xi'(n)=\pi(g)\xi(s(g))$.
Choose a neighborhood $U\subset M$ of $n$ such that $n'\in U$ implies 
$|\la\eta(n'),\xi'(n')\ra_{\HH_{n'}}-\la\eta(n),\xi'(n)\ra_{\HH_{n}}|<\varepsilon/2$. This is possible since $\la\eta,\xi'\ra$ is continuous on $M$. Since $\pi$ is strongly
continuous there exists an open set $V\subset G$ containing $g$ such that for all $g'\in V$ one has $t(g')\in U$
and
\[\|\pi(g')\xi(s(g'))-\xi'(t(g'))\|_{\HH_{t(g')}}<\varepsilon/(2 \sup_{n'\in U}\|\eta(n')\|).\]
Hence, for all $g'\in V$
\[\begin{array}{l}
|\la\eta(t(g')),\pi(g')\xi(s(g'))\ra_{\HH_{t(g')}}-\la\eta(n),\xi'(n)\ra_{\HH_{n}}|\\
\quad\leq|\la\eta(t(g')),\pi(g')\xi(s(g'))\ra_{\HH_{t(g')}}-\la\eta(t(g')),\xi'(t(g'))\ra_{\HH_{t(g')}}|\\
\quad\quad+|\la\eta(t(g')),\xi'(t(g'))\ra_{\HH_{t(g')}}-\la\eta(n),\xi'(n)\ra_{\HH_{n}}|\\
\quad<\|\eta(t(g'))\|\varepsilon/(2  \sup_{n'\in U}\|\eta(n')\|)+\varepsilon/2\leq\varepsilon.
\end{array}\]

The converse implication is proven as follows. Suppose $(\pi,\HH,\Delta)$ is weakly continuous unitary.
Let $U(\varepsilon,\eta,V)$ be a neighborhood of $\pi(g)\xi(s(g))$ in $\HH$ for a given $g\in G$ and $\xi\in\Delta$,
where $\eta\in\Delta$ satisfies $\eta(t(g))=\pi(g)\xi(t(g))$. We compute for any $g'\in G$,
\begin{align}
\|\eta(t(g'))-\pi(g')\xi(s(g'))\|_{\HH_{t(g')}}&=|\la\eta(t(g')),\eta(t(g'))\ra-\la\eta(t(g')),\pi(g')\xi(s(g'))\ra\nonumber\\
&\quad-\la\pi(g')\xi(s(g')),\eta(t(g'))\ra+\la\pi(g')\xi(s(g')),\pi(g')\xi(s(g'))\ra|^{1/2}\nonumber\\
&\leq (|\la\eta(t(g')),\eta(t(g'))\ra-\la\eta(t(g')),\pi(g')\xi(s(g'))\ra|\nonumber\\
&\quad+|\la\xi(s(g')),\xi(s(g'))\ra-\la\pi(g')\xi(s(g')),\eta(t(g'))\ra|)^{1/2} \label{eqn10}
\end{align}
By weak continuity we can choose a neighborhood $W_{g}\subset G$ of $g$ such that $g'\in W_{g}$ implies
\[|\la\eta(t(g')),\pi(g')\xi(s(g'))\ra-\la\eta(t(g)),\pi(g)\xi(s(g))\ra|<\varepsilon.\]
Since $t$ is open and $\eta\in\Delta$, we can choose a $W_{g}'\subset W_{g}$ such that 
\[|\la\eta(t(g')),\eta(t(g'))\ra-\la\eta(t(g)),\eta(t(g))\ra|<\varepsilon\]
Hence the first two terms of Equation (\ref{eqn10}) are smaller than $2\varepsilon$.
Analogously, the last two terms of Equation (\ref{eqn10}) are also smaller than $2\varepsilon$, which finishes the proof.
\end{proof}

\begin{lemma}\label{normalstrong}
If a representation $(\pi,\HH,\Delta)$ is continuous, then it is strongly continuous. The converse holds if $\pi$ is unitary.
\end{lemma}
\begin{proof}
Suppose $(\pi,\HH,\Delta)$ is continuous. Suppose $g\in G$ and $\xi\in\Delta$. Let a neighborhood
$U(\varepsilon,\eta,V)\subset\HH$ of $\pi(g)\xi(s(g))$ be given such that $\eta(t(g)=\pi(g)\xi(s(g))$. Then, by
continuity of $\pi$ there exists a neighborhood $W_{g}\subset G_{s}\times_{p}\HH$  of $g$ such that $g'\in W_{g}$ implies
$\Psi(W_{g})\subset U(\varepsilon,\eta,V)$. Now, define a subset of $G$
\[W_{G}:=\{g'\in G\mid(g',\xi(s(g')))\in W_{g}\}.\]
This set is open since it equals $s^{-1}\xi^{-1}p_{2}(W)\cap p_{1}(W)$. If $g'\in W_{G}$, then
\[\|\eta(t(g'))-\pi(g')\xi(s(g'))\|<\varepsilon.\]

Conversely, suppose $(\pi,\HH,\Delta)$ is strongly continuous and unitary. Suppose $(g,h)\in G_{s}\times_{p}\HH$. Let
$U(\varepsilon,\eta,V)$ be a open neighborhood of $\pi(g)h$ with $\eta(t(g))=\pi(g)h$ as usual. Let $\xi$ be any section
in $\Delta$ such that $\xi(s(g))=h$. Then by strong continuity there exists an open set $V_{g}\subset G$ such that $g'\in V_{g}$ implies $\|\eta(t(g'))-\pi(g')\xi(s(g'))\|<\varepsilon$. Define the set
\[W_{g,h}:=\{(g',h')\in G_{s}\times_{p}\HH\mid \|h'-\xi(s(g'))\|<\varepsilon, g'\in V_{g}\}.\]
It is easily seen to be open and $(g',h')\in W_{g,h}$ implies
\[\begin{array}{ll}
\|\eta(t(g'))-\pi(g')h'\|&\leq\|\eta(t(g'))-\pi(g')\xi(s(g'))\|+\|\pi(g')\xi(s(g'))-\pi(g')h'\|\\
&<\varepsilon+\|\pi(g')\|\|\xi(s(g'))-h'\|<2\varepsilon,
\end{array}\]
which finishes the proof.
\end{proof}

\begin{definition}\label{morphrep}
A \textbf{morphism of continuous (unitary) representations} $(\HH^{1},\Delta^{1},\pi_{1})\to(\HH^{2},\Delta^{2},\pi_{2})$ of a
groupoid is a morphism $\Psi:(\HH^{1},\Delta^{1})\to(\HH^{2},\Delta^{2})$ of
continuous fields of Hilbert spaces that intertwines the groupoid
representations
\[\xymatrix{\HH^{1}_{s(g)}\ar[r]^{\pi_{1}(g)}\ar[d]_{\Psi_{s(g)}}&\HH^{1}_{t(g)}\ar[d]^{\Psi_{t(g)}}\\ 
\HH^{2}_{s(g)}\ar[r]_{\pi_{2}(g)}&\HH^{2}_{t(g)}.}\]
\end{definition}

\begin{example}
The trivial representation of a groupoid $G$ is given by 
the continuous field $(\HH,\Delta)$ that has fiber $\CC$ over each $m\in M$ and a map
$\pi:G\to U(M\times\CC)\simeq M\times U(\CC)\times M$,
\[g\mapsto(t(g),1,s(g)).\]

We give another example of a continuous unitary representation of a groupoid.
For any continuous function $f:G\to\RR$ we can construct the representation
\[\pi_{f}:g\mapsto (t(g),e^{2\pi i(f(t(g))-f(s(g)))},s(g)).\]
These representation are all isomorphic. Indeed, let $f,g:G\to\RR$ then
\[m\mapsto e^{2\pi i (f(m)-g(m))}\]
is an isomorphism $(\HH,\Delta,\pi_{g})\to(\HH,\Delta,\pi_{f})$. In particular all these representations
are isomorphic to $\pi_{0}$, which is
the trivial representation.
\end{example}

Recall that a groupoid is \textbf{proper} if $t\times s: G\to M\times M$ is a proper map.
\begin{proposition}\label{makeun}
If $G$ is a proper groupoid endowed with a Haar system, then any continuous representation 
$(\HH,\Delta,\pi)$ is isomorphic to a unitary representation.  
\end{proposition}
\begin{proof}
Suppose $(\HH,\Delta,\pi)$ is a non-zero continuous representation of $G$.
Let $c:M\to\RR_{>0}$ be a cutoff function (cf. \cite{tuxul}, with $t$ and $s$ interchanged). It exists since $G$ is proper. Define an inner product 
$\la.,.\ra^{new}$ on $\HH$ by the following description: for 
all $m\in M$ and $h,h'\in \HH_{m}$,
\[\la h,h'\ra^{new}(m):=\int_{G_{m}}\la\pi(g)h,\pi(g)h'\ra c(t(g)\lambda_{m}(dg).\]
This inner product is $G$-invariant, since the Haar system and $t$ are right invariant. It gives rise to a new topology of $\HH$. The isomorphism is 
the identity on $\HH$, which is easily seen to be 
continuous. Indeed, let $h\in\HH$ and let $U(\varepsilon,\xi,V)\ni h$ be an open set in $\HH$ with respect to the old norm. Then there exists a  
an open set $V'$ such that $V'\subset V$ and $g\mapsto\|\pi(g)\|$ is bounded on $t^{-1}V'\cap\supp(c\circ t)$. Since $c\circ t$ has compact support, 
the  function 
\[m'\mapsto \int_{g\in G_{m'}}\|\pi(g)\|c(t(g))\lambda_{m'}(dg)\]
is bounded on $V'$. Hence we can set
\[\delta:=\frac{\varepsilon}{\sup_{m\in V'}\int_{g\in G_{m'}}\|\pi(g)\|c(t(g))\lambda_{m'}(dg)}.\]
Then $h'\in U(\delta,\xi,V')$ (in the old topology) implies
\begin{align*}
\|h'-\xi(m')\|^{new}_{m'}&=\int_{G_{m'}}\|\pi(g)(h'-\xi(m')\| c(t(g))\lambda_{m'}(dg)\\
&\leq\int_{G_{m'}}\|\pi(g)\| c(t(g))\lambda_{m'}(dg)\|(h'-\xi(m')\|\\
&\leq\varepsilon,
\end{align*}
which proves the continuity of the identity map.

The proof that the inverse (also the identity) is continuous proceeds similarly. One uses that
\begin{align*}
\|h'-\xi(m')\|&=\int_{G_{m'}}\|h'-\xi(m')\|c(t(g)))\lambda_{m'}(dg)\\
&=\int_{G_{m'}}\|\pi(g^{-1})\pi(g)(h'-\xi(m'))\|c(t(g)))\lambda_{m'}(dg)\\
&=\sup_{g\in G_{m'}}\|\pi(g)\|\int_{G_{m'}}\|\pi(g)(h'-\xi(m'))\|c(t(g)))\lambda_{m'}(dg)\\
\end{align*}
and local boundedness of $g\mapsto\|\pi(g)\|$. This finishes the proof.
\end{proof}

A representation $(\HH,\Delta,\pi)$ is \textbf{locally trivial} if the continuous field $(\HH,\Delta)$ is locally trivial.
In \cite{tuxul} locally trivial representations of a groupoid $G\rra M$ are called $G$-vector bundles.
 
\subsection{Continuity of representations in the operator norm.}\label{opno}
In this section we go through quite some effort to define a suitable topology on the set of bounded linear operators
$\{P:\HH_{m}\to\HH_{n}\}_{n,m\in M}$ for a continuous field of Hilbert spaces $(\{\HH_{m}\}_{m\in M},\Delta_{\HH})$.
This is done not only to be able to consider representations which are continuous in the operator topology, but the continuous
field of Banach spaces thus obtained also plays a crucial r\^{o}le in Section \ref{grconv}.
At first reading one could consider skipping the proofs. 

Suppose $G\rra M$ is a continuous groupoid and let $R_{G}:=(t\times s)(G)\rra M$ be the \textbf{orbit relation groupoid}. It is continuous if
$s,t:G\to M$ are open maps. Let
$(\{\HH_{m}\}_{m\in M},\Delta_{\HH})$ be a continuous field of Hilbert spaces over $M$.
Consider the continuous field of Banach spaces over $R_{G}$ whose fiber at $(n,m)$ is given by the bounded linear operators
$\HH_{m}\to\HH_{n}$, i.e. $\BB_{(n,m)}:=\BB(\HH_{m},\HH_{n})$. This is indeed a Banach space for the norm
\[\|P\|=\sup_{h\in\HH_{m},\|h\|_{\HH_{m}}=1}\|P(h)\|_{\HH_{n}}.\]
We define a space of sections $\Delta_{\BB}$ of the field to be those maps $(n,m)\mapsto P(n,m)$ in $\prod_{(n,m)\in
R}\BB(\HH,\HH)$ such that 
\begin{itemize}
\item[i)] for every $m\in M$ and $h\in\HH_{m}$
\[n\mapsto P(n,m) h\]
is in $\Delta_{\HH}$ , 
\item[ii)] for every $n\in M$ and $\xi\in\Delta_{\HH}$ the map
\[m\mapsto P(n,m)\xi(m)\]
is continuous $M\to\HH_{n}$,
\item[iii)] The map $(n,m)\mapsto\|P(n,m)\|$ locally bounded, and
\item[iv)] $P$ is adjointable, which means that there exists a $P^{*}:R\to \BB(\HH,\HH)$, satisfying i), ii) and iii),  such that for all
$\xi,\eta\in\Delta_{\HH}$ one has $(\eta,P\xi)=(P^{*}\eta,\xi)$, more concretely: for all $(n,m)\in R$
\[\la\eta(n),P(n,m)\xi(m)\ra_{\HH_{n}}=\la P^{*}(m,n)\eta(n),\xi(m)\ra_{\HH_{m}}.\]
\end{itemize}
\begin{lemma}\label{oppietoppie}
The pair $(\{\BB(\HH_{n},\HH_{m})\}_{(n,m)\in R_{G}},\Delta_{\BB})$ is a lower semi-continuous field of Banach spaces.
\end{lemma}
\begin{proof}
We first show that this is true for $R=M\times M$. Then the lemma easily follows since the above field is the restriction of the field to
$R$ (which is a closed subspace of $M\times M$), hence again a continuous field.

First, we prove lower semi-continuity of the norm of a section $P\in\Delta_{\BB}$. This follows from the fact that the map
\[(n,m,h)\mapsto\|P(n,m)h\|_{\HH_{n}}\]
is a continuous map $M\times M\times_{p}\HH\to\RR$, analogously to the proof of Lemma \ref{lsc}. This last statement is proven as follows. Let $\varepsilon>0$ be given.
Suppose $(n,m,h)\in M\times M\times_{p}\HH$. There exists a $\xi\in\Delta_{\HH}$ such that $\xi(m)=h$. Then by condition
(i), (ii), (iii) and continuity of $\|\xi\|$, there exists a neighborhood $W\in M\times M\times_{p}\HH$ such that for any
$(n',m',h')\in W$ the map $\|P\|$ is bounded on $W$ and we have
\[\begin{array}{l}
|\|P(n',m')h'\|-\|P(n,m)h\| |\leq|\,\|P(n',m')h'\|-\|P(n,m')h'\|\,|\\
\quad\quad+|\,\|P(n,m')h'\|-\|P(n,m')\xi(m')\|\,|+|\,\|P(n,m')\xi(m')\|-\|P(n,m)\xi(m)\|\,|\\
\quad\leq\varepsilon+\|P(n,m')\|\varepsilon+\varepsilon.
\end{array}\]

Next, we prove that for every $P\in\BB(\HH_{n},\HH_{m})$ and every $\varepsilon>0$ there exist a $Q\in\Delta_{\BB}$ such
that $\|Q(n,m)-P\|<\varepsilon$. Suppose $P\in\BB(\HH_{n},\HH_{m})$ and let $\varepsilon>0$ be given.
Let $\xi_{1},\ldots,\xi_{k}\in\Delta_{\HH}$ be such that for any $h\in\HH_{m}$
\[\|h-\sum_{i=1}^{k}\la\xi_{i}(m),h\ra\xi_{i}(m)\|<\varepsilon.\]
Let $\eta_{1},\ldots,\eta_{l}\in\Delta_{\HH}$ be such that for any $h\in\HH_{n}$
\[\|h-\sum_{i=1}^{l}\la\eta_{i}(n),h\ra\xi_{i}(n)\|<\varepsilon.\]
Define, for $(n',m')\in M\times M$,
\[Q(n',m')h:=\sum_{i=1}^{k}\sum_{j=1}^{l}\la\xi_{i}(m'),h\ra\la\eta_{j}(n),P \xi_{i}(m)\ra\eta_{j}(n')\]
One easily checks that $Q\in\Delta_{\BB}$.
Furthermore,
\[\begin{array}{l}
\|P h-Q(n,m)h\|\leq \|P h-\sum_{i=1}^{k}\la\xi_{i}(m),h\ra P \xi_{i}(m)\|\\
\quad\quad+\|\sum_{i=1}^{k}\la\xi_{i}(m),h\ra P \xi_{i}(m)-\sum_{i=1}^{k}\sum_{j=1}^{l}\la\xi_{i}(m'),h\ra\la\eta_{j}(n),P\xi_{i}(m)\ra\eta_{j}(n')\|\\
\quad<\|P\|\varepsilon+\varepsilon.
\end{array}\]

The last step is to show that $\Delta_{\BB}$ is locally uniformly closed. Suppose 
\[Q\in\prod_{(n,m)in M\times M}\BB(\HH_{n},\HH_{m}).\] 
Suppose that for all $\varepsilon>0$ and all $(n,m)\in M\times M$ there is a $Q'\in\Delta_{\BB}$
such that
\[\|Q(n',m')-Q'(n',m')\|<\varepsilon\]
on a neighborhood $V$ of $(n,m)$. We shall now show that this implies $Q\in\Delta_{\BB}$. Indeed, let $\varepsilon>0$ be
given and suppose $n\in M$. Then there exist $Q'$ and $V$ as above. Define $U:=p_{1}(V)$. Then $n'\in U$ implies, for
any $h\in\HH_{m}$, that
\[\|Q(n',m)h-Q'(n',m)h\|\leq\|Q(n',m)-Q'(n',m)\|\|h\|<\varepsilon\|h\|.\]
Hence $n\mapsto\|Q(n,m)h\|$ is continuous.
In a similar way one proves condition (ii) for $Q$ which finishes the proof.
\end{proof}

We shall see in Lemma \ref{BBlem} that $\BB(\HH,\HH)$ is a so-called lower semi-continuous Fell bundle over the orbit relation groupoid $R_{G}$ and 
therefore a lower semi-continuous $C^{*}$-category. The collection of sets 
\[\{U(\varepsilon,\xi,V)\mid\xi\in\Delta_{\BB},\varepsilon>0, V\subset R\mbox{ open}\},\] 
as defined in Lemma \ref{toponfield} for a continuous field of Banach spaces, is generally a subbasis for the topology on 
$\coprod_{(n,m)\in R}\BB(\HH_{n},\HH_{m})$, instead of a basis.
Since the field is not continuous in general, we do not have $\Delta=\Gamma_{0}(R_{G},\BB(\HH,\HH))$. 
Consider the restriction of the total space $\BB(\HH,\HH)$ to the unitary operators, i.e. 
\[U^{op}(\HH):=\coprod_{(n,m)\in R}U(\HH_{m},\HH_{n}),\] 
endowed with the subspace topology.
\begin{lemma}\label{oppietoppie2}
The total space $U^{op}(\HH)$ is a continuous groupoid over $M$.
\end{lemma}
\begin{proof}
We show that the composition $\BB(\HH,\HH)^{(2)}\to\BB(\HH,\HH)$ is a continuous map.
First note that for every $(P,Q)\in\BB(\HH,\HH)^{(2)}$ the inequality $\|P Q h\|\leq\|P\|\|Q h\|$ implies
\[\|P Q\|\leq\|P\|\|Q\|.\]
Suppose that $P_{2}\in\BB(\HH_{p},\HH_{n})$, $P_{1}\in\BB(\HH_{n},\HH_{m})$ and $U(\varepsilon,Q,V)$ is an open neighborhood
of $P_{2}P_{1}$ such that $Q(p,m)=P_{2}P_{1}$. There are $Q_{1},Q_{2}\in\Delta_{\BB}$ such that $Q_{1}(n,m)=P_{1}$ and
$Q_{2}(p,n)=P_{2}$. Choose $\varepsilon_{i}>0$ and an open subset $V_{i}\subset M$ such that $P_{i}'\in U(\varepsilon_{i},Q_{i},V_{i})$
implies $\|P_{i}'\|\varepsilon_{i}<\varepsilon/3$ for $i=1,2$. Furthermore, note that by condition (i), for each $m'\in M$ and
$h\in\HH_{m'}$ the map $n'\mapsto Q_{1}(n',m')h$ is in $\Delta_{\BB}$. Hence by condition (ii) the map for each $p',m'\in M$
the map $n'\mapsto Q_{2}(p',n')Q_{1}(n',m')$ is continuous. The map $(p',m')\mapsto Q_{2}(p',n')Q_{1}(n',m')$ is easily seen to be continuous too. Hence we
can shrink $V_{1}$ and $V_{2}$ such that $(p',n',m')\in V_{2}\times_{M} V_{2}$ implies 
\[|\,\|Q_{2}(p',n)Q_{1}(n,m')-Q_{2}(p',n')Q_{1}(n',m')\|\,|<\varepsilon/3.\]
Define $Q\in\Delta_{\BB}$ by $Q(p',m'):=Q_{1}(p',n)Q_{2}(n,m')$
Suppose 
\[(P_{2}',P_{1}')\in U(\varepsilon_{2},Q_{2},V_{2})_{s}\times_{t}U(\varepsilon_{1},Q_{1},V_{1}),\] 
then
\begin{align*}
\|P_{2}'P_{1}'-Q(p',m')\|&=\|P_{2}'P_{1}'-Q_{2}(p',n)Q_{1}(n,m')\|\\
&\leq\|P_{2}'P_{1}'-Q_{2}(p',n')P_{1}'\|+\|Q_{2}(p',n')P_{1}'-Q_{2}(p',n')Q_{1}(n',m')\|\\
&\quad+\|Q_{2}(p',n')Q_{1}(n',m')-Q_{2}(p',n)Q_{1}(n,m')\|\\
&<\|P_{2}'-Q_{2}(p',n')\|\,\|P_{1}'\|+\|Q_{2}(p',n')\|\,\|P_{1}'-Q_{1}(n',m')\|+\varepsilon/3\\
&<\varepsilon_{2}\|P_{1}'\|+\|Q_{2}(p',n')\|\varepsilon_{2}+\varepsilon/3<\varepsilon.
\end{align*}
Proving that the other structure maps are continuous is similar, but easier.
\end{proof}

\begin{definition}
A representation $(\pi,\HH,\Delta)$ is \textbf{continuous in the operator norm} if the map
\[G\to \BB(\HH,\HH),g\mapsto\pi(g)\]
is continuous.
If $G$ is unitary, then the representation is continuous if
\[G\to U^{op}(\HH),g\mapsto\pi(g)\]
is a continuous map of groupoids.
\end{definition}

\begin{lemma}\label{opnormal}
A representation is continuous if it is continuous in the operator norm.
The converse implication is true if the representation $\Delta^{\pi}$ is finitely generated over $C_{0}(M)$ and unitary.
\end{lemma}
\begin{proof}
Suppose $(g,h)\in G_{s}\times_{p}\HH$ and let $n=t(g)$ and $m=s(g)$. Suppose $U(\varepsilon,V,\xi)$ is a neighborhood of $\pi(g)h$, with 
$\xi(n)=\pi(g)h$. Let $Q\in\Delta_{\BB}$ be any section with $Q(n,m)=\pi(g)$, which exists since $(\BB(\HH,\HH),\Delta_{\BB})$ is a lower 
semi-continuous field of Banach spaces. Let $\eta\in\Delta_{\HH}$ be a section such that $\eta(m)=h$. By the conditions i), ii) and iii) above there 
exists a neighborhood $S\subset R$ of $(n,m)$ such that for all $(n',m')\in S$
\[\|\xi(n')-Q(n',m)h\|<\varepsilon/4,\]
the function $\|Q\|$ is bounded on $S$ and
\[\|Q(n',m)\eta(m)-Q(n',m')\eta(m')\|<\varepsilon/4.\]
Define
\[\begin{array}{l}
\delta:=\frac{\varepsilon}{4\sup_{(n',m)\in S}\|Q(n',m')\|},\\
W':=U(\delta,\eta ,p_{2}(S)),\\
K:=\sup_{h'\in W'}\|h'\|,
\end{array}\]
and
\[W:=\pi^{-1}(U(Q,\frac{\varepsilon}{4 K} ,S)),\]
where $p_{2}:M\times M\to M$ is the projection on the second entry.
We claim that $(g',h')\in W_{s}\times_{p}W'$ implies $\pi(g')h'\in U(\varepsilon,V,\xi)$. Indeed, suppose \mbox{$(g',h')\in W_{s}\times_{p}W'$} and $m'=s(g')$, $n'=t(g')$, then
\[\begin{array}{ll}
\|\xi(n')-\pi(g')h'\| & \leq\|\xi(n')-Q(n',m)h\|+\|Q(n',m)\eta(m)-Q(n',m')\eta(m')\| \\
&+\|Q(n',m')\eta(m')-Q(n',m')h'\|+\|Q(n',m')h'-\pi(g')h'\|\\
&<\varepsilon/4+\varepsilon/4+\|Q(n',m')\|\delta+\|h'\|\frac{\varepsilon}{4 K}<\varepsilon.
\end{array}\]

We shall now prove the converse implication. Suppose $(\HH^{\pi},\Delta^{\pi},\pi)$ is a strongly continuous unitary representation on a continuous field of Hilbert spaces with $\Delta^{\pi}$ finitely generated. There exist a finite set $\{\xi_{i}\}_{i\in I}$ of sections in $\Delta^{\pi}$ such that for each $m'\in M$ the set $\{\xi_{i}(m')\}_{i\in I}$ contains a (normalized) basis for $\HH_{m'}$.
Suppose $U(\varepsilon,Q,V)$ is a neighborhood of $\pi(g)$, $s(g)=m$, $t(g)=n$ and $Q(n,m)=\pi(g)$.
Note that by condition (i) $n'\mapsto Q(n',m)\xi_{i}(m)$ is in $\Delta^{\pi}$, so by strong continuity of $\pi$ there exists an open set $U_{i}\subset G$ such that $g'\in U_{i}$ implies
\[\|\pi(g')\xi'(s(g'))-Q(t(g'),m)\xi_{i}(m)\|<\varepsilon/(2 |I|).\]
Moreover, by condition (ii) we can shrink $U_{i}$ such that $g'\in U_{i}$ implies that
\[\|Q(t(g'),m)\xi_{i}(m)-Q(t(g'),s(g'))\xi_{i}(s(g')\|<\varepsilon/(2 |I|).\]
Hence
\[\|\pi(g')\xi'(s(g'))-Q(t(g'),s(g'))\xi_{i}(s(g')<\varepsilon/|I|\]
for $g'\in U_{i}$. Define $U:=\bigcap_{i\in I}U_{i}$, then $g'\in U$ implies
\[\begin{array}{l}
\|\pi(g')-Q(t(g'),s(g'))\|\\
\quad=\sup_{h'\in\HH_{s(g')},\|h'\|=1}\|\pi(g')h'-Q(t(g'),s(g'))h'\|_{\HH_{t(g')}}\\
\quad<\sum_{i\in I}\|\pi(g')\xi_{i}(s(g'))-Q(t(g'),s(g'))\xi_{i}(s(g')\|_{\HH_{t(g')}}\\
\quad<\sum_{i\in I}\varepsilon/|I|=\varepsilon,
\end{array}\]
which finishes the proof.
\end{proof}
From these comparision Lemmas (Lemma \ref{strongweak}, Lemma \ref{normalstrong} and Lemma \ref{opnormal}) we can conclude that 
for unitary representations any of these topologies are equivalent. Hence from now on we shall not specify which notion we mean, but only say that a unitary representation is continuous (if it is).

\subsection{Example: the regular representations of a groupoid.}\label{secregrep}
The following example considers the regular representation. In a different form it was already studied by Reneault (cf.\ \cite{ren}), but this was on 
$L^{2}(G)$ as a measurable field of Hilbert spaces. We are interested in representations on continuous fields. 
Therefore, the statement of Lemma \ref{regrep} is actually new. It generalizes the analogous statement for groups.

Suppose a continuous groupoid $G\rra M$ is endowed with a \textbf{left Haar system}, i.e.\ a left $G$-invariant continuous family of measures
$\{\lambda^{m}\}_{m\in M}$ for $t:G\to M$, cf. Example \ref{exfield}. The left $G$-invariance means that for every $m,n\in M$, $g'\in G_{m}^{n}$ and 
every $f\in C_{c}(G)$
\[\int_{G^{m}}f(g' g)\lambda^{m}(dg)=\int_{G^{n}}f(g)\lambda^{n}(dg).\]
\begin{lemma}\label{regrep}
The \textbf{left regular representation} of $G$ on $(\hat{L^{2}_{t}}(G),\Delta_{t}^{2}(G))$ defined by (continuous extension of)
\[(\pi_{L}(g)f)(g')=f(g^{-1}g'),\]
for $f\in C_{c}(G^{s(g)})$ and $g'\in G^{t(g)}$, is a strongly continuous and unitary representation.
\end{lemma}
\begin{proof}
Unitarity is immediate from the $G$-invariance of the Haar system.

We have to check that for all $\xi\in \Delta_{t}^{2}(G)$ the map
$g\mapsto\pi_{L}(g)\xi(s(g))$ is continuous $G\to\hat{L^{2}_{t}}(G)$. Let $g\in G$.
Suppose a neighborhood $U(\varepsilon,\eta,V)\subset \hat{L^{2}_{t}}(G)$ of $\pi_{L}(g)\xi(s(g))$ is given, where $\varepsilon>0$, $V$
an open set in $M$ and $\eta\in\Delta_{t}(G)$ is a section satisfying $\pi_{L}(g)\xi(s(g))=\eta(s(g))$.
There exist $\xi',\eta'\in C_{c}(G)$ such that $\|\eta-\eta'\|_{\hat{L}^{2}}<\varepsilon/3$, $\|\xi-\xi'\|_{\hat{L}^{2}}<\varepsilon/3$
and $\pi_{L}(g)\xi'(s(g))=\eta'(s(g))$.
To continue we first need the following lemma due to A.\ Connes \cite{con1}.
\begin{lemma}
If $f$ is a compactly supported continuous function on $G^{(2)}$, then the map
\[g\mapsto\int_{h\in G^{s(g)}}f(g,h)\lambda^{s(g)}(d h)\]
is continuous on $G$. 
\end{lemma}
We restate the proof for completeness.
\begin{proof}
Since $G^{(2)}$ is closed in $G\times G$, there exists a continuous and bounded extension $\bar{f}$ of $f$ to $G\times G$ (we suppose here that
$G$ is a normal space). The map $(g,m)\to\int_{h\in G^{m}}f(g,h)\lambda^{m}(d h)$ is continuous, as is proven as follows. Let $(g,m)$ be any
element in $G\times M$ and let $\varepsilon'>0$
be given. Since the Haar system is continuous and $\bar{f}(g',\cdot)$ converges uniformly to
$\bar{f}(g,\cdot)$ for $g'\to g$.
we can choose a neighborhood $W\in G\times M$ such that $(g',m')\in W$ implies
\[\left|\int_{h\in G^{m}}\bar{f}(g',h)\lambda^{m}(d h)-\int_{h'\in G^{m'}}\bar{f}(g',h')\lambda^{m'}(d h')\right|<\varepsilon'/2\]
and
\[\left|\int_{h'\in G^{m'}}\bar{f}(g',h')\lambda^{m}(d h')-\int_{h'\in G^{m'}}\bar{f}(g,h')\lambda^{m'}(d h')\right|<\varepsilon'/2.\]
As a consequence,
\[\begin{array}{l}
\left|\int_{h\in G^{m}}\bar{f}(g',h)\lambda^{m}(d h)-\int_{h'\in G^{m'}}\bar{f}(g,h')\lambda^{m'}(d h')\right|\\
\leq \left|\int_{h\in G^{m}}\bar{f}(g',h)\lambda^{m}(d h)-\int_{h'\in G^{m'}}\bar{f}(g',h')\lambda^{m'}(d h')\right|\\
\quad +\left|\int_{h'\in G^{m'}}\bar{f}(g',h')\lambda^{m}(d h')-\int_{h'\in G^{m'}}\bar{f}(g,h')\lambda^{m'}(d h')\right|\\
< \varepsilon'/2 + \varepsilon'/2=\varepsilon'.
\end{array}\]
Restricting to $\{(g,m)\mid s(g)=m\}\subset G\times M$ gives the required result.
\end{proof}
Now, apply this lemma to the map
\[f(g',h'):=|\xi'((g')^{-1}h')-\eta'(h')|^{2}.\]
As a result, 
\[f(g'):=\sqrt{\int_{h'\in G^{t(g')}}|\xi'((g')^{-1}h')-\eta'(h')|^{2}\lambda^{t(g')}(d h')}\]
depends continuously on $g'$. Note that $f(g)=0$, so that we can choose a neighborhood $U\subset G$ of $g$ such that $f(g')$ is smaller
than $\varepsilon/3$ if $g'\in U$. Finally, intersect $U$ with $t^{-1}(V)$ to obtain the required open set in $G$ whose image is a subset of
$U(\varepsilon,\eta,V)$.
\end{proof}
In the same way one proves that the \textbf{right regular representation} of $G$ on $(\hat{L^{2}_{s}}(G),\Delta_{s}^{2}(G))$ given by
\[\pi_{L}(g)h(g'):=h(g' g)\]
(where $h\in L^{2}(G_{s(g)})$ and $g'\in G_{t(g)}$) is strongly continuous and unitary.

\subsection{Example: continuous families of groups.}\label{famgr}
The following example can give the reader a feeling for the the issues on the global topology with continuous groupoid representations.
We express the set of finite-dimensional continuous representations of a family of groups on a given continuous field of Hilbert spaces in terms of sections of a bundle (or family) of the sets of finite-dimensional continuous representations of each of the groups.

Suppose $H$ is a locally compact group. Let $\RRep(H)$ denote the set of non-zero continuous unitary representations of $H$. 
This set can be endowed with a topology. Indeed, one uses the Jacobson topology on the primitive spectrum of the $C^{*}$-algebra $C^{*}(H)$. 
We shall not go into the details, since there is an easier description of the case that has our interest. Denote by $\RRep^{n}(H)$ the subspace 
of continuous non-zero unitary representations on $\CC^{n}$ with standard inner product $\la z,z'\ra=\bar{z}z'$. 
\begin{lemma}(\cite{dix}, 18.1.9)\label{topfam}
For every integer $n\geq 0$ a subbasis for 
the topology on $\RRep^{n}(H)$ is given by the sets 
\[U(\pi,\varepsilon,K):=\{\pi'\in\RRep(H)\mid\max_{g\in K}|\langle h',\pi(g)\,h\rangle-\langle h',\pi'(g) \,h\rangle|<\varepsilon,\forall 
h,h'\in S(\CC^{n})\},\]
for compact sets $K\subset H$, representations $\pi$ and $\varepsilon>0$. The set $S(\CC^{n})$ is the unit sphere in $\CC^{n}$.
 \end{lemma}

Suppose $s:G\to M$ is a continuous family of groups, i.e.\ $G=\coprod_{m\in M}G_{m}\to M$. 
{\em Fix} a finite-dimensional continuous field of Hilbert spaces $(\HH,\Delta)$. Also choose for each $m\in M$ a group $H_{m}\simeq G_{m}$ and an isomorphism $\psi_{m}:G_{m}:=s^{-1}(m)\to H_{m}$, fixing the group structure at each fiber. Endow $\coprod_{m\in M}H_{m}$ with the topology such that
\[\coprod_{m\in M}\psi_{m}: G\to \coprod_{m\in M}H_{m}\]
is a homeomorphism. Denote the canonical projection $\coprod_{m\in M}H_{m}\to M$ by $s'$.

Suppose $\{U_{i}\}_{i\in I}$ is a covering of $M$.
For any open set $U_{i}$, define
\[\RRep^{\HH}(G|_{U_{i}}):=\coprod_{m\in U_{i}}\RRep^{\dim(\HH_{m})}(H_{m})\]
and the canonical projection
\[p_{i}:\RRep^{\HH}(G|_{U_{i}})\to U_{i}.\]
Suppose 
\[\{\phi_{i}:\HH|_{U_{i}}\hookrightarrow U_{i}\times\CC^{\dim(\HH|_{U_{i}})}\}_{i\in I}\] 
is a local trivialization of $(\HH,\Delta)$, cf.\ Definition \ref{loctriv}. Define for all $i,j\in I$ the homeomorphism
\[\gamma_{ij}:=\phi_{j}(\phi_{i})^{-1}:\im(\phi_{i})|_{U_{i}\cap U_{j}}\to\im(\phi_{j})|_{U_{i}\cap U_{j}}.\]

We need the following technical notion. Suppose $p:N\to M$ is a continuous map.
We say a set $K\subset N$ is \textbf{$p$-open-compact}\index{open-compact sets for a map} if the restriction $K\cap p^{-1}(m)$ is compact for all 
$m\in M$ and
the image $p(K)\subset M$ is open.
We say that $p:N\to M$ is \textbf{locally open-compact}\index{locally open-compact map} if every $n\in N$ has a $p$-open-compact neighborhood.
\begin{example}
If $p:N\to M$ is a fiber bundle with locally compact fiber, then it is easy to show that $p$ is locally open-compact. 
\end{example}

For each $i\in I$ the following sets form a subbasis of a topology on $\RRep^{\HH}(G)|_{U_{i}}$:
For any $\xi,\eta\in\Delta$, $V\subset\CC$ open and $K\subset\coprod_{m\in U_{i}}H_{m}$ $s'$-open-compact, 
\[U(\xi,\eta,K,V):=\left\{\pi\in \RRep^{\HH}(G|_{U_{i}})\mid \la\xi,\pi\eta\ra(K\cap H_{p_{i}(\pi)})\subset V\right\}.\]
Define
\[\RRep^{\HH}(G):=(\coprod_{i\in I}\RRep^{\HH}(G|_{U_{i}}))/\sim,\]
where $\RRep^{\HH}(G|_{U_{i}})|_{U_{i}\cap U_{j}}\ni\pi_{i}\sim\pi_{j}\in\RRep^{\HH}(G|_{U_{j}})|_{U_{i}\cap U_{j}}$ iff $\pi_{j}= \gamma_{ij}\pi_{i}\gamma_{ij}^{-1}$.
The space $\RRep^{\HH}(G)$ is uniquely determined up to homeomorphism by the chosen covering $\{U_{i}\}_{i\in I}$ and isomorphisms $\{\psi_{m}:G_{m}\to H_{m}\}_{m\in M}$.

One easily sees that $s:G\to M$ being locally open-compact implies that the topology of $\RRep^{\HH}(G|_{U_{i}})$ restricted to each fiber is equivalent to the topology of Lemma \ref{topfam}.
\begin{proposition}
Suppose that $s:G\to M$ is locally open-compact family of groups.
Then there is a one-to-one correspondence between continuous representations of $s:G\to M$ on $(\HH,\Delta)$ and continuous sections of $\RRep^{\HH}(G)$.
\end{proposition}
\begin{proof}
A continuous unitary representation $\pi$ of $G$ on $(\HH,\Delta)$ corresponds to a continuous section of $\RRep^{\HH}(G)$, i.e.\
to a family of sections $\tilde{\pi}_{i}:U_{i}\to\RRep^{\HH}(G)|_{U_{i}}$ given by 
\[\tilde{\pi}_{i}(m)=\phi_{i}\circ\pi\circ(\psi_{m}^{-1}\times\phi_{i}^{-1}).\]
These are easily seen to be compatible, i.e. $\tilde{\pi}_{j}= \gamma_{ij}\tilde{\pi}_{i}\gamma_{ij}^{-1}$. It remains to show that each
$\tilde{\pi}_{i}$ is continuous. Consider an open set $U(\xi,\eta,K,V)$ as above. Note that
\begin{eqnarray*}
\tilde{\pi}^{-1}(U(\xi,\eta,K,V))&=&\{m\in U_{i}\mid \la\xi,\pi\eta\ra|_{K\cap H_{m}}\subset V\}\\
&=&s'(K\cap \{g\in \coprod_{m\in U_{i}}H_{m}\mid \la\xi,\pi\eta\ra\subset V\}),
\end{eqnarray*}
which is open since $K$ is $s'$-compact and $\pi$ is continuous.

A continuous section $\tilde{\pi}$ of $\RRep^{\HH}(G)$ determines a continuous unitary representation by
\[\pi(g):=\phi_{i}^{-1}\circ\tilde{\pi}_{i}\circ(\psi_{s(g)}(g)\times\phi_{i})\in U(\HH_{s(g)}),\]
where $i\in I$ such that $s(g)\in U_{i}$.
We only need to show that $\pi|_{G_{U_{i}}}\psi^{-1}$ is continuous. Suppose $\xi,\eta\in\Gamma^{0}(\im(\phi_{i}))$. Given $g\in \coprod_{m\in U_{i}}H_{m}$ and $V\subset\CC$, let
$K$ be an $s'$-compact neighborhood of $g$ and $W\subset K$ an open neighborhood of $g$. Consider $U(\xi,\eta,K,V)$. 
Define $W':=W\cap s^{-1}\tilde{\pi}^{-1}(U(\xi,\eta,K,U))$, which is open since $s$ and $\tilde{\pi}$ are
continuous. Then $g'\in W'$ implies
\[\la\xi(s(g')),\pi(g')\eta(s(g'))\ra=\la\xi(s(g')),\tilde{\pi}(s(g'))(g')\eta(s(g'))\ra\in U,\]
which finishes the proof.
\end{proof}

\begin{example}
Consider a locally compact group $H$ and a continuous principal $H$-bundle $\tau:P\to M$. From this we can construct a continuous bundle of groups $P\times_{H}H\to M$,
where the action of $H$ on $H$ is given by conjugation.
Consider a local trivialization \mbox{$\{\chi_{i}:P|_{U_{i}}\to U_{i}\times H\}_{i\in I}$} of $P\to M$. Suppose $I=\NN$. One can fix the group structure at each fiber
of $P\times_{H}H\to M$ as follows: for every $m\in M$ choose the smallest $i\in I$ such that $m\in U_{i}$ and define 
\[\psi_{m}:(P\times_{H}H)_{m}\to H, [p,h]\mapsto \chi_{i}(p)\,h\,\chi_{i}(p)^{-1}.\]
Given a representation $(\pi,\CC^{n})\in\RRep^{n}(H)$, one can construct a vector bundle $\HH:=P\times_{\pi}\CC^{n}\to M$.
Obviously, the trivialization of $P\to M$ gives rise to a trivialization\linebreak\mbox{$\{\phi_{i}:\HH|_{U_{i}}\to U_{i}\times\CC^{n}\}_{i\in I}$} of $\HH\to M$, by
$\phi_{i}([p,z])=(\tau(p),\pi(\chi_{i}(p))z)$. Using these data one can form the bundle $\RRep^{\HH}(P\times_{H}H)\to M$ and a topology on it. 
A continuous section of this bundle is given by
\[\tilde{\pi}_{i}(m)=(h\mapsto\pi(\gamma_{i j}^{-1}h\gamma_{i j}),\]
for all $i\in\NN$, $m\in U_{i}$, $h\in H$ and the smallest $j\in\NN$ such that $m\in U_{j}$.
This section corresponds to the representation of $P\times_{H}H$ given by $\tilde{\pi}([p,h])[p,z]=[p,\pi(h)z]$. 
\end{example}

\begin{remark}\label{twist}
One can ``twist'' $\HH:=P\times_{\pi}\CC^{n}$ by another continuous field $(\HH',\Delta')$, carrying the trivial representation of 
$P\times_{H}H\to M$, to obtain a representation on $\HH\otimes\HH'$. A similar construction is possible for any groupoid, cf.\ Lemma \ref{modstruc}.
\end{remark}

\subsection{Representations of the global bisections group.}
For the reader who prefers representation theory of groups and wonders why one should be interested in representations of groupoids at all, the next section will be of particular interest. Namely, to any continuous groupoid is associated a topological group: the group of global bisections.
For a large class of continuous groupoids (the ones we call locally bisectional) we establish a bijection between continuous representation of the groupoid on continuoous fields of Hilbert spaces and a specific type of continuous representations of the group of global bisections on Banach spaces.
Hence the representation theory of such groupoids can be ``embedded'' in the representation theory of groups. From this point of view, the groupoid offers a way to study the some representations of these groups of bisections.

Suppose $G\rra M$ is a continuous groupoid. A \textbf{global bisection} is a map $\sigma:M\to G$ such that $t\circ\sigma=id_{M}$ and
$\tilde{\sigma}:=s\circ\sigma:M\to M$ is a homeomorphism. Denote the set of global bisections of $G$ by $\Bis(G)$. This set has a group
structure, cf.\ \cite{weca}. Moreoever, it is even a topological group. 
\begin{lemma}
$\Bis(G)$ has the structure of a topological group in the compact-open topology.
\end{lemma}
\begin{proof}
The multiplication is given by
\[(\sigma_{1}\cdot\sigma_{2})(m):=\sigma_{1}(m)\sigma_{2}(\tilde{\sigma}_{1}(m)).\]
The unit is given by the unit section $u:M\to G$ and the inverse is defined by
\[\sigma^{-1}(m):=(\sigma(\tilde{\sigma}^{-1}(m)))^{-1}.\]
The group laws are easily checked, for example
\[\begin{array}{ll}
(\sigma\cdot\sigma^{-1})(m)&=\sigma(m)\sigma^{-1}(\tilde{\sigma}(m))\\
&=\sigma(m)(\sigma(\tilde{\sigma}^{-1}\tilde{\sigma}(m)))^{-1}\\
&=1_{m}.
\end{array}\]
We prove that multiplication is continuous $\Bis(G)\times\Bis(G)\to\Bis(G)$. Suppose\linebreak $\sigma_{1}\cdot\sigma_{2}\in U(C,V)$, where $C$ is a
compact set in $M$, $V$ open in $G$ and $U(C,V)$ the set of maps \mbox{$\tau:M\to G$} that satisfy $\tau(C)\subset V$, i.e. $U(C,V)$ is in the standard
sub-basis of the topology on $\Bis(G)$. For each $m\in C$, let $V_{m}$ be a neighborhood of
$(\sigma_{1}\cdot\sigma_{2})(m)=\sigma_{1}(m)\sigma_{2}(\tilde{\sigma_{1}}(m))$. These $V_{m}$ cover $\sigma_{1}\cdot\sigma_{2}(C)$ which is compact
by continuity of the multiplication in $G$ and $\sigma_{1},\sigma_{2}$. Let $\{V_{i}\}_{i\in I}$ be a finite sub-cover. The inverse image $m^{-1}(V_{i})$
is open and contains a Cartesian product of opens $W_{i}^{1}\times W_{i}^{2}$ for each $i\in I$. Then $\sigma_{1}'\in U(C,\bigcup_{i\in
I}W_{i}^{1})$ and $\sigma_{2}'\in U(\tilde{\sigma_{1}}(C),\bigcup_{i\in I}W_{i}^{2})$ implies $\sigma_{1}'\cdot\sigma_{2}'\in U(C,V)$.
\end{proof}
\begin{example} 
The global bisection group of the pair groupoid $M\times M$ is the group of homeomorphisms of $M$. 
\end{example}
\begin{example}
For any group
bundle $G\times M$ over $M$ the group of global bisections is just the group of sections with the pointwise multiplication. In particular, if $M$
is the circle $\mathbb{S}^{1}$ and $G$ a Lie group then the group of global bisections is the loop group $C(\mathbb{S}^{1},G)$ with its usual topology 
(cf.\ \cite{prse}).
\end{example}
\begin{lemma}
A continuous unitary representation $(\pi,\HH,\Delta)$ of a groupoid $G\rra M$ canonically induces a continuous isometric representation of $\Bis(G)$ 
on $\Delta$.
\end{lemma}
\begin{proof}
Define the representation $\tilde{\pi}$ of $\Bis(G)$ by
\[(\tilde{\pi}(\sigma)\xi)(m):=\pi(\sigma(m))\xi(\tilde{\sigma}(m)),\]
where $\xi\in\Delta$, $m\in M$ and $\sigma\in\Bis(G)$.
This representation is isometric, since $\pi$ is unitary:
\[\|\tilde{\pi}(\sigma)\xi\|=\sup_{m\in M}\|\pi(\sigma(m))\xi(\tilde{\sigma}(m)))\|_{\HH_{m}}=\|\xi\|.\]
Continuity is proven as follows. Suppose $\varepsilon>0$ and $\xi\in\Delta$ are given. There exists a compactly supported section
$\xi'\in\Delta_{c}:=C_{c}(M)\Delta$ such that $\|\xi-\xi'\|<\varepsilon/6$. Denote the support of $\xi'$ by $K$.
Moreover, since $\pi$ is continuous and unitary it is norm continuous and hence there exists an open set $V\subset G$ such 
that $g,g'\in V$ implies $\|\pi(g)\xi'(s(g))-\pi(g')\xi'(s(g'))\|<\varepsilon/3$. Now, suppose that $\sigma,\sigma'\in U(K,V)$ and
$\eta\in B(\xi,\varepsilon/6)$, then
\[\sup_{m\in M}\|\pi(\sigma(m))\eta(\tilde{\sigma}(m))-\pi(\sigma'(m))\eta(\tilde{\sigma}'(m))\|<\varepsilon,\]
which finishes the proof.
\end{proof}
The obtained representation of $\Bis(G)$ is actually \textbf{$C_{0}(M)$-unitary}, in the sense that
\[\la\tilde{\pi}(\sigma)\xi,\tilde{\pi}(\sigma)\eta\ra=\la\xi,\eta\ra\]
for all $\sigma\in\Bis(G)$ and $\xi,\eta\in\Delta$.

For the following result we need a technical condition on groupoids.
We call a continuous groupoid $G\rra M$ \textbf{bisectional} if
\begin{itemize}
\item[(1)] every $g\in G$ is in the image of a bisection;
\item[(2)] for all compact sets $K\subset M$ and open sets $V\subset G$, the set $\bigcup_{\sigma\in U(K,V)}\im(\sigma)\subset G$ is open.
\end{itemize}

\begin{theorem}\label{bisect}
Suppose $G\rra M$ is bisectional. Then
there is a bijective correspondence between continuous unitary representations of $G$ and continuous $C_{0}(M)$-unitary representations of $\Bis(G)$ on a Hilbert $C_{0}(M)$-algebra satisfying 
\begin{itemize}
\item $C_{0}(M)$-linearity, i.e.
\[\tilde{\pi}(\sigma)(f\,\xi)=\tilde{\sigma}^{*}f\,\tilde{\pi}(\sigma)(\xi)\]
for all $\sigma\in\Bis(G)$, $\xi\in\Delta$ and $f\in C_{0}(M)$
and 
\item locality, i.e. if $\sigma(m)=1_{m}$ for some $m\in M$, then $\|\tilde{\pi}(\sigma)\xi-\xi\|(m)=0$
\end{itemize}
\end{theorem}
\begin{proof}
Given a representation $(\tilde{\pi},\Delta)$ of $\Bis(G)$ as above, define a representation \mbox{$\pi:G\to U(\HH)$} as follows.
Form the continuous field of Hilbert spaces $\{\HH_{m}\}_{m\in M}$ associated to $\Delta$. For any $g\in G$ and $h\in\HH_{s(g)}$,
define
\[\pi(g)h:=(\tilde{\pi}(\sigma)\xi)(t(g)),\]
for any $\xi\in\Delta$ such that $\xi(s(g))=h$ and $\sigma\in\Bis(G)$ such that $\sigma(t(g))=g$, which exist by assumption.
We now show that this definition does not depend on the choice of $\sigma$ and $\xi$. Suppose $\xi,\xi'$ satisfy $\xi(m)=h=\xi'(m)$.
Let $\{U_{i}\}_{i\in\NN}$ be a family of sets such that $\bigcap_{i\in\NN}U_{i}=\{s(g)\}$ and $\{\chi_{i}:U_{i}\to [0,1]\}$ a family of functions
such that $\chi_{i}(s(g))=0$ and $\chi_{i}(n)=1$ for all $n\in M\backslash U_{i}$. Then
\begin{eqnarray*}
(\tilde{\pi}(\sigma)\xi)(t(g))-(\tilde{\pi}(\sigma)\xi')(t(g))&=&\lim_{i\to\infty}(\tilde{\pi}(\sigma)\chi_{i}(\xi-\xi'))(t(g))\\
&=&\lim_{i\to\infty}\chi_{i}(\tilde{\sigma}(t(g)))(\tilde{\pi}(\sigma)(\xi-\xi'))(t(g))\\
&=&0,
\end{eqnarray*}
since $\tilde{\pi}$ is $C_{0}(M)$-linear and $\tilde{\sigma}(t(g))=s(g)$.

Suppose $\sigma(m)=\sigma'(m)$ for $\sigma,\sigma'\in\Bis(G)$ and $m\in M$. Then, by locality, for all $\xi\in\Delta$
\[\|\tilde{\pi}(\sigma^{-1}\,\sigma')\xi-\xi\|(m)=0,\]
and hence $(\tilde{\pi}(\sigma)\xi)(m)=(\tilde{\pi}(\sigma')\xi)$.

Unitarity of $\pi$ follows at once from $C_{0}(M)$-unitarity of $\tilde{\pi}$.

Next, we prove continuity of $\pi$. Suppose $(g,h)\in G\fiber{s}{p}\HH$ and $U(\varepsilon,\eta,V)$ open neighborhood of $\pi(g)h$, where
$\eta(t(g))=\pi(g)h$. We need to construct an open neigborhood of $(g,h)$, which maps to $U(\varepsilon,\eta,V)$.
Consider
\[B(\eta,\varepsilon):=\xi\in\Delta\mid \|\eta-\xi\|<\varepsilon.\]
Let a $\sigma\in\Bis(G)$ be such that $\sigma(t(g))=g$, which exists since $G$ is bisectional.
Define $\xi:=\bar{\pi}(\sigma)^{-1}\eta$. By continuity of $\bar{\pi}$ there exists an open neighborhood $B(\xi,\delta)$ of $\xi$ and
an open neighborhood $U(K,W)$ of $\sigma$ such that $\bar{\pi}(U(K,W)\times B(\xi,\delta))\subset B(\eta,\varepsilon)$.
Since $G\rra M$ is bisectional, there exists an open neighborhood $W'$ of $g$ in $\bigcup_{\sigma\in U(K,W)}\im(\sigma)$.

Suppose that $(g',h')\in W'\fiber{s}{p}U(\xi,\delta,\tilde{\sigma}^{-1}(V))$, then
\[\pi(g')h'=(\bar{\pi}(\sigma')\xi')(t(g'))\in U(\varepsilon,\eta,V),\]
for some $\sigma'\in U(K,W)$ and $\xi'\in B(\xi,\delta)$.

One easily sees that the constructions given in this proof to obtain representations of $G$ from representations of $\Bis(G)$ and vice versa
in the proof of the above lemma are inverses of each other. 
\end{proof}

\section{Groupoid representation theory}
Is there a Schur's Lemma for groupoids? Is there a Peter-Weyl theorem for groupoids? In this section we give answers to these questions. We discuss a way to generalize these statements to groupoids. It turns out that you need extra conditions on the groupoid for 
the statements to be true (unlike what is suggested in \cite{amin}). A crucial 
r\^{o}le is played by the functors that restict representations of a groupoid to representations of its isotropy groups. This section shows that 
representation theory of groupoids is quite different from representation theory for groups, but many results can be carried over using some caution.

\subsection{Decomposability and reducibility}
\begin{definition}
The direct sum of a family of continuous representations $\{(\HH^{i},\Delta_{i},\pi_{i})\}_{i\in I}$
of a groupoid $G\rra M$ is defined as follows. The family of Hilbert spaces is given by  $\HH_{m}:=\bigoplus_{i\in
I}\HH^{i}_{m}$. The space $\Delta$ is the smallest Banach space containing all finite sums of sections $\sum_{j\in
J}\xi_{j}$, where $\xi_{j}\in\Delta_{j}$, such that $(\HH,\Delta)$ is a continuous field of Hilbert spaces.
The representation of $G$ on $(\HH,\Delta)$ is given by extending the map 
$\bigoplus_{j\in J}\pi_{j}:g\mapsto\sum_{j\in J}\pi_{j}(g)$ on finite sums.
\end{definition}
We say that that a continuous unitary representation $(\HH,\pi)$ of a groupoid
$G$ is \textbf{decomposable} if it is equivariantly isomorphic to
a direct sum of representations $(\HH^{1},\pi_{1})$ and $(\HH^{2},\pi_{2})$
\[\HH\simeq\HH^{1}\oplus\HH^{2}.\]
and \textbf{indecomposable} if this is not possible. 

A \textbf{continuous subfield} of a continuous field of Hilbert spaces $(\HH,\Delta)$ is a continuous field
$(\HH',\Delta')$, such that $\HH'_{m}\subset\HH_{m}$ is a closed linear subspace with the induced inner product for all
$m\in M$ and $\Delta'\subset\Delta$. 
\begin{definition}
A \textbf{continuous subrepresentation} of a continuous unitary representation
$(\HH,\pi)$ of a groupoid $G$ is a continuous subfield of $(\HH,\Delta)$ stable under
$\pi$. 
\end{definition}
\begin{proposition}
If $(\HH,\Delta,\pi)$ is a continuous \textit{locally trivial} unitary representation and \linebreak\mbox{$(\HH',\Delta',\pi')$} a 
\textit{locally trivial} subrepresentation of $(\HH,\Delta,\pi)$, then $(\HH,\Delta,\pi)$ decomposes as a direct sum of $(\HH',\Delta',\pi')$
and another subrepresentation. 
\end{proposition}
\begin{proof}
For each $m\in M$ let $\HH_{m}''$ be the orthogonal complement with respect to the inner product. The family $\{\HH_{m}''\}_{m\in M}$ forms a continuous
field, with 
\[\Delta'':=\{\xi\in\Delta\mid \xi(m)\in\HH_{m}''\mbox{ for all }m\in M\},\]
since $\HH$ is locally trivial.
Moreover, $(\HH'',\Delta'')$ is locally trivial too. 
Since $\pi$ is unitary, this complement is $G$-invariant. 
\end{proof}

A continuous unitary representation is \textbf{reducible} if it has a proper continuous subre\-pre\-sen\-tation.
It is \textbf{irreducible} if it is not reducible. 
\begin{example}\label{exdis}
A simple example is a family of groups over a discrete set. A representation of such a family is irreducible iff it has support on one point, 
where it is an irreducible representation of the group at that point.
\end{example}
Decomposability implies reducibility (irreducible implies indecomposable), but not vice versa. Indeed, a representation can contain a 
subrepresentation without being decomposable. 
\begin{example}\label{subfd}
For example, consider the trivial representation of $\RR\rra\RR$ on $(\RR\times\CC,C_{0}(\RR))$. 
It has a subrepresentation given by the continuous field of Hilbert
spaces which is $0$ at $0$ and $\CC$ elsewhere, with space of sections
\[C^{0}_{0}(M):=\{f\in C_{0}(M)\mid f(0)=0\}.\]
This subrepresentation has no complement, since this would be a field that is
$\CC$ at $0$ and zero elsewhere, whose only continuous section could be
the zero section. 
Note that $\RR\rra\RR$ is an example of a groupoid which has no continuous irreducible representations.
\end{example}
Define the \textbf{support of a continuous field of Hilbert spaces} $(\HH,\Delta_{\HH})$ by
\[\supp(\HH,\Delta_{\HH}):=\{m\in M\mid\HH_{m}\not=0\}.\]
This last set equals
\[\{m\in M\mid\xi(m)\not=0\mbox{ for some }\xi\in\Delta_{\HH}\}.\]
One easily sees that for all continuous fields of Hilbert spaces $(\HH,\Delta_{\HH})$ the support $\supp(\HH,\Delta_{\HH})$ is open in $M$.
\begin{lemma}
If the support of a continuous unitary representation $(\HH,\Delta_{\HH},\pi)$ of a groupoid $G$ properly contains a closed union of orbits, 
then it is reducible.
\end{lemma}
\begin{proof}
Let $(\HH,\Delta_{\HH},\pi)$ be a continuous representation of $G$. Suppose $G m\subset M$ is a closed orbit.
Define a new continuous field of Hilbert spaces by
\[\HH_{m}':=\left\{\begin{array}{ll}\HH_{m}&\mbox{ if }m\notin G m\\0&\mbox{ if }m\in G m\end{array}\right.\]
and
\[\Delta_{\HH'}:=\{\xi\in\Delta\mid\xi|_{G m}=0\},\]
The groupoid $G$ represents on $(\HH',\Delta')$ by 
\[\pi'(g):=\left\{\begin{array}{ll}\pi(g)&\mbox{ if }s(g)\notin G m\\id_{0}&\mbox{ if }m\in G m\end{array}\right.\]
One easily sees that $(\HH',\Delta_{\HH'},\pi')$ is a continuous subrepresentation of $(\HH,\Delta_{\HH},\pi)$.
\end{proof}
The representation $(\HH',\Delta_{\HH'},\pi')$ is called the restriction of $(\HH,\Delta_{\HH},\pi)$ to $(G m)^{c}$.
\begin{example}
If a groupoid $G$ is proper (i.e.\ $t\times s:G\to M\times M$ is proper), then its orbits are closed. Hence an irreducible
representation must consist of one orbit, which is clopen, since it is the support of a continuous field and the orbit of a proper groupoid. 
Therefore, a Hausdorff space $M\rra M$ has an irreducible representation iff it has a discrete point $m\in M$.
\end{example}
Along the same lines one can easily show:
\begin{lemma}
If the support of a representation properly contains a clopen set closed under $G$, then the representation is decomposable.
\end{lemma}
As an example, consider Example \ref{exdis}.

\subsection{Schur's lemma.}\label{secsl}
In the previous section we have seen that in many cases of interest the irreducible representations do exist. Therefore, we introduce the weaker 
notion of $M$-irreducibility. A continuous representation $(\pi,\HH,\Delta)$ of a groupoid $G\rra M$ is called \textbf{$M$-irreducible} if the 
restriction of 
$\pi$ to each of the isotropy groups is an irreducible representation. Obviously, if a representation is irreducible, then it is $M$-irreducible. The 
converse does not hold as we have seen in Example \ref{subfd}.
\begin{example}\label{exgrbu}
Suppose $H$ is a topological group, $P\to M$ a continuous principal $H$-bundle and $(\pi,V)$ an irreducible representation of $H$. Then, 
$P\times_{H}V\to M$ carries a canonical $M$-irreducible (but reducible) representation of the bundle of groups 
$P\times_{H}H\to M$ (cf.\ Section \ref{famgr}).
\end{example}
\begin{example}
If $M$ is a topological space with a non-trivial rank 2 vector bundle $E\to M$. 
Then $E\to M$ is not $M$-irreducible as a representation of $M\rra M$, even though it might be indecomposable.
\end{example}
\begin{example}A morphism of $M$-irreducible continuous representations is not necessarily an
isomorphism or the zero map, even if the restriction to each isotropy group is an irreducible representation.
A counterexample is given by the following: let $G$ be the constant bundle of
groups $\RR\times U(1)\rra\RR$. It
represents $M$-irreducible on the trivial rank one vector bundle $\HH:=\RR\times\CC$
over $\RR$ by scalar multiplication. The map $\Psi:(x,z)\mapsto
(x,x\cdot z)$ is an equivariant adjointable map $\HH\to\HH$, not equal to a
scalar times the identity or zero. 
\end{example}
What one does see in this example is that $\Psi$ is a function times the
identity on $\HH$, namely the function $\lambda:\RR\to\CC, x\mapsto x$,
i.e. $\psi=\lambda 1_{\HH}$. An alternative formulation of Schur's lemma for
groupoids would be that an endomorphism of an $M$-irreducible
representation $(\HH,\pi)$ is a function $\lambda\in C(M)$ times the identity on
$\HH$.

For a continuous groupoid $G\rra M$ denote the set of isomorphism classes of continuous unitary representations by $\RRep(G)$. 
Denote the subset of isomorphism classes of indecomposable representation by $\IdRep(G)$, the subset isomorphism classes of of irreducible representations by $\IrRep(G)$ and the set of isomorphism classes of $M$-irreducible representations by $M$-$\IrRep(G)$.
\begin{lemma}[Schur's Lemma for groupoids]
Suppose $(\pi_{i},\HH^{i},\Delta^{i})$ is an $M$-irreducible representation for $i=1,2$.
\begin{itemize}
\item[i)] every equivariant endomorphism $\Psi:\HH^{1}\to\HH^{1}$ is equal to a continuous function $\lambda\in
C(M)$ times the identity on $E$, i.e. $\psi=\lambda 1_{\HH^{1}}$.
\item[ii)] If $\Phi:\HH^{1}\to\HH^{2}$ is a morphism of representations then $\Phi_{m}$ is either an
isomorphism or the zero map $\HH^{1}_{m}\to\HH^{2}_{m}$ for all $m\in M$.
\item[iii)] If, furthermore, $\Res_{m}:M$-$\IrRep(G)\to\IrRep(G_{m}^{m})$ is injective for every $m\in M$, then 
\[\Hom_{G}(\HH^{1},\HH^{2})=\left\{\begin{array}{ll}\mbox{a line bundle}&\mbox{ if }(\pi_{1},\HH^{1},\Delta^{1})\simeq(\pi_{2},\HH^{2},\Delta^{2});\\
0&\mbox{ if }(\pi_{1},\HH^{1},\Delta^{1})\not\simeq(\pi_{2},\HH^{2},\Delta^{2}).\end{array}\right.\]
\end{itemize}
\end{lemma}
The proof follows easyly from the analogous statement for groups.
\begin{example}
Suppose $P\to M$ is a principal $H$-bundle for a group $H$. If $G\rra M$ is the gauge groupoid $P\times_{H}P\rra M$,
then every irreducible representation is $M$-irreducible. Moreover, $\Res_{m}:M$-$\IrRep(G)\to\IrRep(G_{m}^{m})$ is injective for all $m\in M$.
\end{example}
\begin{example}
Consider the two-sphere as a groupoid $S^{2}\rra S^{2}$. It is proper and all indecomposable vector bundles over $S^{2}$ have rank one. These are $M$-irreducible representations, but obviously $\Res_{m}:M$-$\IrRep(S^{2})\to\IrRep(\{m\})$ is not injective for any $m\in M$.
\end{example}

\begin{corollary}\label{mcoef}
If a continuous groupoid $G$ has the property that for all $m\in M$ the restriction map
\[\Res_{m}:M\mbox{-}\IrRep(G)\to\IrRep(G_{m}^{m})\]
is injective,
then for any two non-isomorphic $M$-irreducible unitary representations $(\HH,\Delta,\pi)$, $(\HH',\Delta',\pi')$
and \mbox{$\xi,\eta\in\Delta$}, \mbox{$\xi',\eta'\in\Delta'$},
\[\la\la\xi,\pi\eta\ra,\la\xi',\pi'\eta'\ra\ra_{\hat{L}^{2}(G)}=0\]
\end{corollary}
\begin{proof}
This easily follows from the version of this statement for compact groups and the invariance of the Haar system.
\end{proof}

\subsection{Square-integrable representations.}\label{square}
In this section we define the notion of square-integrability for continuous groupoid representations. In the end, we prove that for proper groupoids,
with $M/G$ compact, unitary representations are square-integrable, generalizing an analogous result for compact groups.

Suppose $G\rra M$ is a locally compact groupoid endowed with a Haar system
$\{\lambda_{m}\}_{m\in M}$, which desintegrates as
$\lambda_{m}=\int_{n\in t(G_{m})}\lambda_{m}^{n}\mu_{m}(dn)$, for a Haar system $\{\mu_{m}\}_{m\in M}$ on $R_{G}\rra M$ and
a continuous family of measures $\{\lambda_{m}^{n}\}_{(n,m)\in R_{G}}$ on $t\times s:G\to M\times M$. 
Using the family $\{\lambda_{m}^{n}\}_{(n,m)\in R_{G}}$ one can construct the continuous field of Hilbert spaces
\[(\hat{L}^{2}(G),\Delta^{2}(G)):=(\hat{L}^{2}_{t\times s}(G),\Delta_{t\times s}^{2}(G)),\] 
over $R_{G}$, cf.\ Example \ref{exfield}.
\begin{example}
Very simple example of this is the following. If $M$ is a space and $\mu$ a measure on it and $H$ a Lie group with Haar measure $\lambda$. Then the trivial transitive groupoid $M\times G\times  M\rra M$ with isotropy groups $H$ has a Haar system 
$\{\lambda_{m}=\lambda \times \mu\}_{m\in M}$. Obviously, this decomposes as $\lambda_{m}=\int_{n\in M}\lambda\mu(dn)$, hence
\[(\hat{L}^{2}(G),\Delta^{2}(G))= (L^{2}(G,\lambda)\times (M\times M), C_{0}(M\times M,L^{2}(G,\lambda)).\]
\end{example}
A map $f:G\to\CC$ is called \textbf{$\hat{L}^{2}(G)$-square integrable} if the induced map
\[(m,n)\mapsto (g\mapsto f(g),G_{m}^{n}\to\CC)\]
is in $\Delta^{2}(G)$.

The \textbf{conjugate representation} $(\bar{\HH},\bar{\Delta},\bar{\pi})$ of a representation $(\HH,\Delta,\pi)$ is defined as follows. The family of Hilbert spaces
is given by $\bar{\HH}_{m}=\HH_{m}$ as Abelian groups, but with conjugate complex scalar multiplication. Also, the space
of sections $\bar{\Delta}=\Delta$ remains the same (but with conjugate $C_{0}(M)$-action). The representation of $G$ on 
$(\Bar{\HH},\Delta)$ is given by $\bar{\pi}(g)h=\pi(g)h$, where $h\in\bar{\HH}_{s(g)}$.

The \textbf{tensor product $(\HH^{1}\otimes\HH^{2},\Delta^{\otimes},\pi_{1}\otimes\pi_{2})$ of two continuous representations}
$(\HH^{1},\Delta^{1},\pi_{1})$ and $(\HH^{2},\Delta^{2},\pi_{2})$
of a groupoid $G$ is defined as follows. The family of Hilbert spaces is given by 
$\HH_{m}:=\HH^{1}_{m}\otimes\HH^{2}_{m}$. The space $\Delta^{\otimes}$ is the smallest Banach space containing all finite sums of
sections $\sum_{j\in J}\xi_{j}\otimes \eta_{j}$ of $\xi_{j}\in\Delta^{1}$ and $\eta_{j}\in\Delta^{2}$, such that
$(\HH,\Delta)$ is a continuous field of Hilbert spaces.
The representation of $G$ on $(\HH,\Delta)$ is given by linearly extending the map 
$(\pi_{1}\otimes\pi_{2})(g)(h\otimes h')=\pi(g)h\otimes\pi(g)h'$.

\begin{definition}
A continuous representation $(\pi,\HH,\Delta)$ is \textbf{square-integrable} if the map
\[(\bar{\HH}\otimes\HH,\Delta^{\otimes})\to(\hat{L^{2}}(G),\Delta^{2}(G))\]
given by
\[h_{2}\otimes h_{1}\mapsto(g\mapsto(h_{2},\pi(g)h_{1})_{\HH_{t(g)}}\]
is a map of continuous fields of Hilbert spaces.
\end{definition}
This means that the \textbf{matix coefficients} $\la\xi,\pi\eta\ra$, defined by
\[(n,m)\mapsto(g\mapsto\la\xi(n),\pi(g)\eta(m)\ra)\]
for $\xi,\eta\in\Delta$ are $\hat{L}^{2}(G)$-square-integrable maps.
\begin{example}
For example, consider a topological space $M$.
A (finite-dimensional) vector bundle $E\to M$ is a square-integrable representation of $M\rra M$.
\end{example}
\begin{example}
Consider the family of continuous groups $G:=(\RR\times\ZZ/2\ZZ)\backslash (0,-1)\rra\RR$.
One easily sees that the trivial representation $g\mapsto id_{\CC}$ on $(\RR\times\CC,C_{0}(\RR))$ is
not square-integrable. But, note that $G$ is not proper (although for every $m\in M$ the 
set \mbox{$s^{-1}(m)=t^{-1}(m)$} is compact).
\end{example}
\begin{lemma}
If $G\rra M$ is proper and $M/G$ compact, then every unitary representation is square-integrable.
\end{lemma}
\begin{proof}
Suppose $(\HH,\Delta,\pi)$ is a unitary representation and $\xi,\eta\in\Delta$. Given $\varepsilon>0$,
choose $\xi',\eta'\in C_{c}(M)\Delta$ such that $\|\xi-\xi'\|<\varepsilon'$ and
$\|\eta-\eta'\|<\varepsilon'$, where
\[\varepsilon'=\frac{\min\{\varepsilon,1\}}{3 M\max\{\|\xi\|,\|\eta\|\}}\]
and 
\[M=\max_{(n,m)\in R_{G}}\lambda_{m}^{n}(G_{m}^{n}),\]
which exists since $M/G$ is compact.
First note that $\la\xi',\pi\eta'\ra$ has compact support, since $G\rra M$ is proper.
Moreover,
\begin{eqnarray*}
\|\la\xi,\pi\,\eta\ra-\la\xi',\pi\eta'\ra\|_{\hat{L}^{2}}&\leq&\|\la(\xi-\xi'),\pi\,\eta\ra\|
+\|\la \xi',\pi\,(\eta-\eta')\ra\|\\
&\leq&\max_{(n,m)\in R_{G}}\lambda_{m}^{n}(G_{m}^{n})(\|\xi-\xi'\|\|\eta\|+\|\xi'\|\|\eta-\eta'\|)\\
&\leq&\varepsilon'\|\eta\|+(\|\xi\|+\varepsilon')\varepsilon'\leq\varepsilon,
\end{eqnarray*}
which finishes the proof.
\end{proof}

\subsection{The Peter-Weyl theorem I.}
Suppose $G\rra M$ is a continuous groupoid endowed with a Haar system $\{\lambda_{m}\}_{m\in M}$, which decomposes using a continuous family of
measure $\{\lambda_{m}^{n}\}_{(n,m)\in R_{G}}$ as in Section \ref{square}. Let $\mathcal{E}(G)\subset\Delta^{2}(G)$ denote the
$C_{0}(R_{G})$-submodule spanned by the matrix coefficients (cf.\ Section \ref{square}) of all finite-dimensional representations of
$G\rra M$.

A generalization of the Peter-Weyl theorem as we are going to prove (cf.\ Theorem \ref{PW1} and Theorem \ref{decomp}) appears not to be true for all
continuous groupoids. Therefore, we introduce an extra condition: 
\begin{definition}
For a continuous groupoid $G\rra M$ the restriction map
\[\Res_{m}:\RRep(G)\to\RRep(G_{m}^{m})\]
is dominant if for every $m\in M$ and every continuous unitary representation $(\pi,V)$ of $G_{m}^{m}$ there exists a continuous unitary representation 
$(\pi',\HH,\Delta)$ of $G$ such that $(\pi,V)$ is isomorphic to a subrepresentation of $(\pi'|_{G^{m}_{m}},\HH_{m})$.
\end{definition}
\begin{example}
Suppose $H$ is a group and $P\to M$ a principal $H$-bundle. Since $(P\times_{H}P)_{m}^{m}\simeq H$ and $P\times_{H}P\rra M$ are Morita equivalent, 
$\Res_{m}:\RRep(P\times_{H}P)\to\RRep((P\times_{H}P)_{m}^{m})$ is dominant for all $m\in M$.
\end{example}
\begin{example}
Suppose $H$ is a compact connected Lie group that acts on manifold $M$. Consider the action groupoid $G:=H\ltimes M\rra M$.
\begin{proposition}\label{dom}
The restriction map $\Res_{m}:\RRep(H\ltimes M)\to\RRep((H\ltimes M)_{m}^{m})$ is dominant for all $m\in M$.
\end{proposition}
\begin{proof}
First we note that from every representation $(\pi,V)\in\RRep(H)$ we can construct a representation $\tilde{\pi}:H\ltimes M\to U(M\times V)$
of $H\ltimes M\rra M$ on $M\times V\to M$ by $\tilde{\pi}(h,m):(m,v)\mapsto (h\cdot m,\pi(h)v)$.
Note that the isotropy groups of $H\ltimes M\rra M$ coincide with the isotropy groups of the action. These are subgroups of $H$, hence the question
is whether every representation of a subgroup of $H$ occurs as the subrepresentation of the restiction of a representation of $H$.

Suppose $K$ is a compact Lie subgroup of $H$. Fix a maximal tori $T_{K}\subset K$ and $T_{H}\subset H$ such that $T_{K}\subset T_{H}$, with Lie algebras $\mathfrak{t}_{K}$ and $\mathfrak{t}_{H}$. Note that
$T_{K}\simeq \mathfrak{t}_{K}/\Lambda_{K}$ and $T_{H}\simeq \mathfrak{t}_{H}/\Lambda_{H}$ for lattices $\Lambda_{K}\subset \mathfrak{t}_{K}$ and
$\Lambda_{H}\subset \mathfrak{t}_{H}$. There is an injective linear map
$M:\mathfrak{t}_{K}\to \mathfrak{t}_{H}$ that induces the inclusion $\mathfrak{t}_{K}/\Lambda_{K}\hookrightarrow \mathfrak{t}_{H}/\Lambda_{H}$.
Let $P_{K}$ denote the integral weight lattice of $T_{K}$ and $P_{H}$ the integral weight lattice of $T_{H}$.
Hence $q:=M^{T}:\mathfrak{t}_{H}^{*}\to \mathfrak{t}_{K}^{*}$ is surjective map, mapping $P_{H}$ onto $P_{K}$. 
Hence restriction of representations $\RRep(T_{H})\to\RRep(T_{K})$ is surjective too, since for tori irreducible 
representations correspond to integral weights.

The following argument is valid if one fixes positive root systems $R_{K}^{+}$, $R_{H}^{+}$ and hence fundamental Weyl chambers $C_{K}^{+}$, 
$C_{H}^{+}$ in a way specified in \cite{heck}.
Suppose $(\pi_{\lambda},V)$ is a an irreducible representation of $K$ corresponding to the dominant weight $\lambda\in P_{K}\cap C_{K}^{+}$. 
One can choose any integral weight $\Lambda\in q^{-1}(\lambda)\cap P_{H}\cap C_{H}^{+}$; this set is non-empty, since $q$ is surjective and the 
positive root systems have been fixed appropriately. Let $\pi_{\Lambda}$ denote the irreducible representation of $H$ associated to 
$\Lambda$.
Then the multiplicity of $\pi_{\lambda}$ in $\pi_{\Lambda}|_{K}$ is a positive integer(not necessarily 1), as follows from the Multiplicity Formula
(3.5) in \cite{heck}. This finishes the proof.
\end{proof}
\end{example}
\begin{example}
A simple, but non-Hausdorff example of a proper groupoid which has a non-dominant restriction map is defined as follows. Consider
$\RR\times\ZZ/2\ZZ\rra\RR$ and identify $(x,0)$ with $(x,1)$ for all $x\not=0$. Endow the obtained family of groups 
$(\RR\times\ZZ/2\ZZ)/\sim\rra\RR$, 
with the quotient topology. The non-trivial irreducible representation of $\ZZ/2\ZZ$ is not in the image of $\Res_{0}:\RRep(G):\to \RRep(\ZZ/2\ZZ)$.
\end{example}

We now prove a generalization of the Peter-Weyl theorem for groupoids. Consider the continuous field of Hilbert spaces $(\hat{L}^{2}(G),\Delta^{2}(G))$ associated to a groupoid $G\rra M$.
Let $\overline{\mathcal{E}(G)}$ denote the closure of $\mathcal{E}(G)$ to a Hilbert $C_{0}(R_{G})$-module.
\begin{theorem}[Peter-Weyl for groupoids I]\label{PW1}
If $G\rra M$ is a proper groupoid, $M/G$ is compact and $\Res_{m}$ is dominant for all $m\in M$, then
\[\overline{\mathcal{E}(G)}=\Delta^{2}(G).\]
\end{theorem}
\begin{proof}
Note that $G_{m}^{m}$ is compact so Peter-Weyl for compact groups applies. Using
the dominance property 
\[\overline{\{\Theta(m,m)|\Theta\in\mathcal{E}(G)\}}=L^{2}(G_{m}^{m},\lambda_{m}^{m}),\]
since $(\HH,\Delta,\pi)<(\HH',\Delta',\pi')$, implies $\la\xi,\pi'\eta\ra=\la\xi,\pi\eta\ra$ for $\xi,\eta\in\Delta$.

Note that $l_{g}^{*}:L^{2}(G_{m}^{m},\lambda_{m}^{m})\to L^{2}(G_{m}^{n},\lambda_{m}^{n})$ is an isometry for a chosen $g\in G_{n}^{m}$. Thus 
$\overline{\{l_{g}^{*}(\Theta(m,m))|\Theta\in\mathcal{E}(G)\}}=L^{2}(G_{m}^{n}).$
But, for all $h\in G_{m}^{n}$ and every continuous unitary finite-dimensional representation $(\HH,\Delta ,\pi)$
\[\begin{array}{ll}
l_{g}^{*}(\xi,\pi\eta)(h)&=(\xi(t(g)),\pi(g h)\eta(s(h)))_{\HH_{t(g)}^{\pi}}\\
&=\sum_{k=1}^{\dim(\HH_{n})}(\xi(m),\pi(g)e_{k}(n))_{E_{n}}(e_{k}(n),\pi(h)\eta(m))_{E_{m}},
\end{array}\]
where $e_{1},\ldots,e_{\dim(\HH_{n})}$ are sections which form a basis of $\HH$ at $n$.
Thus $l^{*}_{g}(\xi,\pi\eta)$ is a linear combination 
of matrix coefficients $(e_{k},\pi\eta)$ restricted to $G_{m}^{n}$, which implies 
$\overline{\{\Theta(n,m)|\Theta\in\mathcal{E}(G)\}}=L^{2}(G_{m}^{n})$.

Let $f\in\Delta^{2}(G)$ and $\varepsilon>0$ be given, then there
exists a section $\tilde{f}\in\Delta^{2}(G)$ with compact support
$K$ such that
$\|f-\tilde{f}\|<\varepsilon/2$, where the norm is the one associated to the
$C_{0}(M)$-valued inner product. Moreover, for all $(m,n)\in R$
there are representations $(\HH_{m,n},\Delta_{m,n},\pi_{m,n})$ and sections
$u_{m,n},v_{m,n}\in\Delta_{m,n}$, such that
\[\|\tilde{f}-(u_{m,n},\pi_{m,n}v_{m,n})\|_{L^{2}(G_{n}^{m})}<\varepsilon/2.\]
Since $\pi_{m,n},u_{m,n}$ and $v_{m,n}$ are continuous we can find an open
neighborhood $S_{m,n}\subset R$, such that still
\[\|\tilde{f}-(u_{m,n},\pi_{m,n}v_{m,n})\|_{\hat{L}^{2}(G)|_{S_{m,n}}}<\varepsilon/2,\]
for all $(m,n)\in R$. These $S_{m,n}$ cover $K$, thus there is a finite
subcover, which we denote by $\{S_{i}\}_{i\in I}$ to reduce the indices.
Denote the corresponding representations by $\pi_{i}$ and sections by $u_{i}$
and $v_{i}$ for $i\in I$. Let $\{\lambda_{i}\}$ be a partition of
unity subordinate to $\{S_{i}\}$. Define $\tilde{u}_{i}=\sqrt{\lambda_{i}}u_{i}$
and $\tilde{v}_{i}=\sqrt{\lambda_{i}}v_{i}$, then 
\[g=\sum_{i\in I}(\tilde{u_{i}},\pi_{i}\tilde{v_{i}})\]
is a finite sum of matrix coefficients and
\[\begin{array}{ll}
\|f-g\|&\leq\|f-\tilde{f}\|+\|\tilde{f}-g\|\\
&\leq\varepsilon/2+\sup_{(m,n)\in R}\|\tilde{f}-\sum_{i\in
I}(\tilde{u_{i}},\pi_{i}\tilde{v_{i}})\|_{L^{2}(G_{n}^{m})}\\
&=\varepsilon/2+\sup_{(m,n)\in R}\|\sum_{i\in I}\lambda_{i}\tilde{f}-\sum_{i\in
I}(\sqrt{\lambda_{i}}u_{i},\pi_{i}\sqrt{\lambda_{i}}v_{i})\|_{L^{2}(G_{n}^{m})}\\
&\leq\varepsilon/2+\sum_{i\in I}\lambda_{i}\sup_{(m,n)\in R}\|\tilde{f}-\sum_{i\in
I}(u_{i},\pi_{i}v_{i})\|_{L^{2}(G_{n}^{m})}\\
&\leq\varepsilon/2+\sum_{i\in I}\lambda_{i}\varepsilon/2=\varepsilon,
\end{array}\]
which finishes the proof.
\end{proof}

\begin{example}
For a space $M$, $\overline{\mathcal{E}(M\rra M)}=C_{0}(M)$ and $\overline{\mathcal{E}(M\times M\rra M)}=C_{0}(M\times M)$ as Theorem \ref{PW1} asserts.
\end{example}
\begin{example}
If $H$ is a compact group and $P\to M$ an $H$-principal bundle. Then, for the bundle of groups $P\times_{H}H\to M$ one finds
(cf.\ Example \ref{exgrbu}),
\begin{align*}
\overline{\mathcal{E}(P\times_{H}H\rra M)}&\simeq \Gamma_{0}(P\times_{H}\overline{\mathcal{E}(H)})\\
&\simeq \Gamma_{0}(P\times_{H}L^{2}(H))\\
&\simeq \Delta^{2}(P\times_{H}H),
\end{align*}
where in the second line we used the Peter-Weyl theorem for the group $H$.
\end{example}

\subsection{The Peter-Weyl theorem II.}
In this section we shall try to find a decomposition analogous to the case of compact groups $H$, where one has 
$L^{2}(H)\simeq\bigoplus_{(\pi,V)\in\hat{H}}\bar{V}\otimes V$ equivariantly. 
There is a seemingly relevant proposition that asserts that 
\begin{proposition}\label{tuxulprop}(\cite{tuxul}, Proposition 5.25)
Any locally trivial representation $(\HH,\Delta,\pi)$ of a proper groupoid $G\rra M$ is a direct summand of the regular representation,
after stabilizing, i.e.\ $\HH\subset\hat{L}^{2}_{s}(G)\otimes\mathbb{H}$, $G$-equivariantly, where $\mathbb{H}$ denotes a standard separable Hilbert space, say $l^{2}(\NN)$. 
\end{proposition}
\begin{example}
The Serre-Swan theorem for vector bundles is a nice example of this. Consider the groupoid $M\rra M$ for a space $M$. 
Locally trivial representations of this groupoid are vector bundles. The theorem states that any vector bundle is a direct summand
of $\hat{L}^{2}(M)\otimes \mathbb{H}\simeq M\times \mathbb{H}$. The Serre-Swan Theorem is actually somewhat stronger, since instead of $\mathbb{H}$ 
one could put a finite-dimensional vector space $\CC^{N}$ for large enough $N\in\NN$.
\end{example}
In general the direct summands will not add up to the whole of $\hat{L}^{2}_{s}(G)\otimes\mathbb{H}$, as one sees in the following Example \ref
{ex3}. Moreover, stabilization is not something that occurs in the case of compact groups (which we want to generalize). Therefore, we have to 
choose a different approach.

The first problem is which $\hat{L}^{2}$ continuous field related to the groupoid one has to use.
\begin{example}\label{ex3}
Consider the pair groupoid $M\times M\rra M$ for a space $M$. It has just one irreducible and indecomposable representation, namely the trivial one
$M\times\CC\to M$. Suppose $\mu$ is a Borel measure on $M$.
Note that $\hat{L}^{2}_{s}(M\times M)\simeq L^{2}(M)\times M$ is hence in general to big. The field $\hat{L}^{2}(M\times M)$ has base space 
$M\times M$ and is not well suited as well, since representations are on continuous fields over $M$.
Recall that the isotropy groupoid of $G\rra M$ is denoted by $I(G)\rra M$. It is a continuous family of groups in the subspace topology.
Now consider the continuous field of Hilbert spaces $\hat{L}^{2}(I(M\times M))\simeq \CC\times M$. This seems to be a good candidate.
\end{example}

The continuous field of Hilbert spaces $(\hat{L^{2}}(I(G)),\Delta^{2}(I(G)))$ is the restriction of
$(\hat{L^{2}}(G),\Delta^{2}(G))$ over \mbox{$R_{G}\subset M\times M$} to the diagonal. It carries a continuous unitary representation
\[\pi_{LR}(g)f(h):=f(g^{-1} h g),\]
where $g\in G_{m}^{n}$, $h\in G_{n}^{n}$ and $f\in L^{2}(G_{m}^{m})$. Remember that this field may not exist, since
suitable measures $\lambda_{m}^{m}$ may not exist.

\begin{lemma}
For any square-integrable continuous unitary representation $(\HH^{\pi},\Delta^{\pi},\pi)$ there is an
equivariant map
\[\Psi_{\pi}:(\bar{\HH^{\pi}}\otimes \HH^{\pi},\Delta^{\otimes})\to(\hat{L}^{2}(I(G)),\Delta^{2}(I(G))),\]
given by
\[h_{2}\otimes h_{1}\mapsto(g\mapsto(h_{2},\pi(g)h_{1})_{\HH_{t(g)}}).\]
\end{lemma}
This map is a slight adaptation of the one introduced for the definition of square-integrability.
\begin{proof}
For equivariance we compute
\[\begin{split}
\Psi(\pi(g)(h_{1}\otimes h_{2}))&=\Psi(\pi(g)h_{1}\otimes\pi(g)h_{2})\\
&=(g'\mapsto (\pi(g)h_{1},\pi(g')\pi(g)h_{2}))\\
&=(g'\mapsto (h_{1},\pi(g^{-1})\pi(g')\pi(g)h_{2}))\\
&=(g'\mapsto (h_{1},\pi(g^{-1} g' g)h_{2}))\\
&=\pi_{LR}(g)(g'\mapsto(h_{1},\pi(g')h_{2}))
\end{split}\]
which finishes the proof.
\end{proof}

\begin{theorem}[Peter-Weyl for groupoids II]\label{decomp}
Suppose $G$ is a proper groupoid and for every $m\in M$
\[\Res_{m}:M\mbox{-}\IrRep(G)\to \IrRep(G_{m}^{m})\]
is bijective. Then
\begin{equation}\label{decomp2}
\bigoplus_{\pi\in M\mbox{-}\IrRep(G)}\Psi_{\pi}:\bigoplus_{\pi\in\hat{G}}(\bar{\HH^{\pi}},\bar{\Delta^{\pi}})\otimes(\HH^{\pi},\Delta^{\pi})\to(\hat{L}^{2}(I(G)),\Delta^{2}(I(G)))
\end{equation}
is an isomorphism of representations.
\end{theorem}
\begin{proof}
The above lemma gives a $G$-equivariant map.
Surjectivity of this map follows from Theorem \ref{PW1}. Injectivity follows from
Corollary \ref{mcoef}.
\end{proof}
\begin{example}
Consider a principal $H$-bundle $P\to M$ for a compact group $H$ and the associated gauge groupoid $G:=P\times_{H}P\rra M$. By Morita equivalence 
(cf.\ Section \ref{morequ}) there is a bijection between unitary irreps $(V,\pi)$ of $H$ and unitary indecomposable, irreducible representations
$P\times_{H}V\to M$ of $G$. Therefore, $\Res_{m}$ is bijective. Hence, by Theorem \ref{decomp}, one has
the decomposition of formula \ref{decomp2}. This is no surprise, since $I(P\times_{H}P)\simeq P\times_{H}H$, where $H$ acts on $H$ by conjugation,
hence
\begin{align*}
\hat{L}^{2}(I(P\times_{H}P)&\simeq P\times_{H}L^{2}(H)\\
&\simeq P\times_{H}\bigoplus_{(\pi,V)\in\hat{H}}\bar{V^{\pi}}\otimes V^{\pi}\\
&\simeq \bigoplus_{(\pi,V^{\pi})\in\hat{H}}(P\times_{H}\bar{V^{\pi}})\otimes(P\times_{H} V^{\pi})\\
&\simeq\bigoplus_{(\pi,\HH^{\pi})\in M\mbox{-}\IrRep(P\times_{H}P)}\bar{\HH^{\pi}}\otimes \HH^{\pi}.
\end{align*}
This is exactly the statement of Theorem \ref{decomp}.
\end{example}
\begin{remark}
Only for a few (types of) groupoids the map $\Res_{m}:M\mbox{-}\IrRep(G)\to \IrRep(G_{m}^{m})$ is bijective for all $m\in M$.
If the map is just surjective, then one could try to find a subset $\PW(G)$ of $M\mbox{-}\IrRep(G)$ that does map bijectively
to $\IrRep(G_{m}^{m})$ for every $m\in M$. We call such a set a \textbf{PW-set}. Then, Theorem \ref{decomp} holds with $M$-$\IrRep(G)$ replaced by
$\PW(G)$.
\end{remark}
\begin{example}
If $H$ is a compact group and $P\to M$ a principle $H$-bundle, then $G:=P\times_{H}H\to M$ is a bundle of groups and
\[\PW(G):=\{P\times_{H}V\mid (\pi,V)\in\IrRep(H)\}\]
is a PW-set (cf.\ Example \ref{exgrbu}). Twisting the representations in $\PW(G)$ with a fixed non-trivial line bundle over $M$ gives another
PW-set, hence these sets are not unique.
\end{example}

\subsection{Morita equivalence.}\label{morequ}
It is well-known that a Morita equivalence of groupoids induces an equivalence of the categories of continuous represenations \textit{on vector 
bundles} of these groupoids. In this section we generalize this to representations on continuous fields of Hilbert spaces.

We begin by brushing up on generalized morphisms and Morita equivalences of continuous groupoids (cf.\ \cite{landsman1, mrw, hisk, mr, mrcun}).
Suppose $G\rra G_{0}$ is a continuous groupoid. Suppose $G$ acts continuously from the left on a map $J:N\to G_{0}$. 
The action is called \textbf{left principal} if the map
\[(g,n)\mapsto(g\cdot n, n)\]
is a homeomorphism
\[G\times_{G_{0}}N\to N\times_{G\backslash N} N,\]
where $G\backslash N$ is endowed with quotient topology.
Suppose $H$ is another continuous groupoid over $H_{0}$. A space $N$ is a \textbf{$G-H$-bibundle} if it carries a left $G$
action and a right $H$ action which commute, i.e. $(g\cdot m)\cdot h=g\cdot(m\cdot h)$, $J_{H}(g\cdot m)=J_{H}(m)$ and $J_{G}(m\cdot
h)=J_{G}(m)$. 
A \textbf{morphism of $G-H$-bibundles} $N$, $N'$ is a $G-H$-equivariant continuous map $N\to N'$.
An isomorphism class of a left principal $G-H$-bibundle can be seen as an arrow $G\to H$ in a 
category of groupoids. The arrows are called \textbf{Hilsum-Skandalis maps} or \textbf{generalized morphisms}. The unit morphism $U(G)$ at $G$ is given by the class of $G$ itself as a $G-G$-bibundle, with left and right multiplication as actions.
One can show that a morphism given by a class of bibundles is an isomorphism if the representing bundles have principal left and
right actions. In that case, one easily sees that $N/H\simeq G_{0}$ and $G\backslash N\simeq H_{0}$. Groupoids which are isomorphic in
this category are called \textbf{Morita equivalent}.
One can prove that a $G-H$ bibundle $N$ represents a Morita equivalence if it is left and right principal and $N/H\simeq G_{0}$ and 
$G\backslash N\simeq H_{0}$.

For a continuous groupoid $G\rra G_{0}$, denote the category of continuous unitary representations of $G\rra G_{0}$ on continuous fields of Hilbert 
spaces by $\Rep(G)$. Before we continue we need a technical tool. Suppose $G$ is a continuous groupoid endowed with a Haar system
$\{\lambda_{m}\}_{m\in M}$.
\begin{definition}
A left action of $G$ on a map $J:N\to G_{0}$ has \textbf{left Dirac sequences} if for each $n\in N$ there exists a sequence of functions
$(\delta_{k}^{n})_{k\in\NN}$ on $N$ such that
\begin{itemize}
\item[i)] $\delta_{k}^{n}\geq 0$ on $J^{-1}(J(n))$,
\item[ii)] $\int_{g\in G_{J(n)}}\delta_{k}^{n}(g\cdot n)\lambda_{n}(d g)=1$ for all $k\in\NN$,
\item[iii)] For every open neighborhood $U\subset G^{m}$ of $1_{m}$ and every $\varepsilon>0$ there is an $N\in\NN$ such that for $k>N$
\[\int_{g\in U^{c}}\delta_{k}^{n}(g\cdot n)\lambda^{n}(d g)<\varepsilon.\]
\end{itemize}
\end{definition}
There is an analogous notion of right Dirac sequences for right actions of groupoids.
 
\begin{theorem}\label{morrep}
Suppose $G$ and $H$ are continuous groupoids endowed with a right Haar system.
If $G$ and $H$ are Morita equivalent and there exists a Morita
$G$-$H$-bibundle that admits left and right Dirac sequences, then the representation categories $\Rep(G)$ and $\Rep(H)$ are equivalent.
\end{theorem}
\begin{proof}
Let $[N]:G\to H$ be a Hilsum-Skandalis map. It gives a map $\RRep(G)\to\RRep(H)$ on the objects as follows.
Let $\pi:G \to U(\HH)$ be a representation of $G$ on a continuous field $(\{\HH_{m}\}_{m\in M},\Delta_{\HH})$ of Hilbert
spaces. Consider the pullback continuous field of Hilbert spaces $J_{G}^{*}(\{\HH_{m}\}_{m\in M},\Delta_{\HH})$ over $N$. 
Note that the projection $J_{G}^{*}(\HH)\to N$ is equivariant. Define 
\[(J^{*}_{G}\Delta_{\HH})_{G}:=\{\xi\in\Gamma_{0}(G\backslash(J^{*}_{G}\HH))\},\]
where $G\backslash(J^{*}_{G}\HH)$ is endowed with the quotient topology.
\begin{lemma}
The pair $(G\backslash(J_{G}^{*}\{\HH_{m}\}_{m\in M}),(J_{G}^{*}\Delta_{\HH})_{G})$ is a continuous field of Hilbert spaces over $G\backslash N\simeq H_{0}$.
\end{lemma}
\begin{proof}
Write $\alpha:G\times_{M}N\to N\times_{G\backslash N} N$ for the homeomorphism $(g,n)\mapsto(g\cdot n,n)$. In particular, there is homeomorphism
$\alpha_{n}:G_{J_{G}(n)}\to G n$. for every $n\in N$. Suppose $[h_{1},n], [h_{2},n]\in (G\backslash(J^{*}\HH))_{G n}$. As a result of the previous remark we
can indeed find unique representatives $h_{1}, h_{2}\in\HH_{n}$ and define
\[\begin{array}{rl}
[h_{1},n]+[h_{2},n]&:=[h_{1}+h_{2},n];\\
\lambda [h_{1},n]&:=[\lambda h_{1},n];\\
\langle [h_{1},n], [h_{2},n]\rangle_{(G\backslash(J^{*}\HH))_{G n}}&:= \langle h_{1},h_{2}\rangle_{\HH_{n}}.
\end{array}\]
Note that $[h,g\cdot n]=[\pi(g)h,n]$. Hence the inner product is well-defined, since $\pi$ is unitary. 
Therefore, every fiber $(G\backslash(J^{*}\HH))_{G n}$ is a Hilbert space.
In fact, one easily sees $(G\backslash(J^{*}\HH))_{G n}\simeq\HH_{n}$.

What is left to prove is the fact that for every $[h,n]\in(G\backslash(J^{*}\HH))_{G n}$ there is a
$\xi\in(J^{*}\Delta_{\HH})_{G}$, such that 
$\xi([h, n])=G h$. This follows from the fact that every $h\in J^{*}\HH$ can be approximated by a sections
$\{\xi'_{k}\}_{k\in\NN}$ in $(J^{*}\Delta_{\HH})_{G}$ which are the image of 
\textit{equivariant} sections $\{\xi_{k}\}_{k\in\NN}$ in $C(N,J^{*}\HH)$ under the projection $J^{*}\HH\to G\backslash J^{*}\HH$.

The construction of the $\xi_{k}$ is as follows. Suppose $h\in\HH_{n}$, then there is an $\eta\in J^{*}\Delta$ such that
$\eta(n)=h$. Moreover, we may suppose that
the support of $\eta$ is compact, by multiplying with a function on $N$ with compact support, which is 1 at $n$.
Let $(\delta^{n}_{k})_{k\in\NN}$ be a Dirac sequence for $N$ at $n$. Define the averaged section by
\[\xi_{k}(n'):=\int_{g\in G^{J(n')}}\delta_{k}^{n}(g^{-1}\cdot n')\pi(g)\eta(g^{-1} n')\lambda^{J(n)}(d g).\]
This integral exists, since the support of $\eta$ is compact and the action of $G$ on $N$ principal.
It depends continuously on $n'$, since the integrand is compactly supported and continuous. Indeed, $\pi$ is strongly
continuous, which implies
that it is weakly continuous. By the properties of the Haar system this (Bochner) integral depends continuously on $n'$.
Moreover, it is equivariant, as follows from the invariance of the Haar system.
Let $\varepsilon>0$ be given. Let $U\subset G^{J(n)}$ be a neighborhood of $1_{J(n)}$, such that for all $g\in U$
\[\|\pi(g)\eta(g^{-1} n)-\eta(n)\|<\varepsilon.\]
Then, there exists a $k\in\NN$ such that 
\[\int_{U^{c}}\delta_{k}^{n}(g\cdot n)\lambda^{J(n)}(d g)<\varepsilon.\]
Hence,
\[\begin{split}
\|\xi_{k}(n)-h\|&=\|\int_{g\in G^{J(n)}}\delta_{k}^{n}(g^{-1} n)\pi(g)\eta(g^{-1} n)\lambda^{J(n)}(d g)-h\|_{(J^{*}\HH)_{n}}\\
&=\|\int_{g\in G^{J(n)}}\delta_{k}^{n}(g)(\pi(g)\eta(g^{-1} n)-\eta(n))\lambda^{J(n)}(d g)\|_{(J^{*}\HH)_{n}}\\
&\leq\int_{g\in G^{J(n)}}\delta_{k}^{n}(g)\|\pi(g)\eta(g^{-1} n)-\eta(n)\|_{(J^{*}\HH)_{n}}\lambda^{J(n)}(d g)\\
&=\int_{g\in U}\delta_{k}^{n}(g)\|\pi(g)\eta(g^{-1} n)-\eta(n)\|_{(J^{*}\HH)_{n}}\lambda^{J(n)}(d g)\\
&\quad+\int_{g\in U^{c}}\delta_{k}^{n}(g)\|\pi(g)\eta(g^{-1} n)-\eta(n)\|_{(J^{*}\HH)_{n}}\lambda^{J(n)}(d g)\\
&\leq \int_{g\in U}\delta_{k}^{n}(g)\lambda^{J(n)}(d g)\varepsilon\\
&\quad+\varepsilon\cdot\max_{g\in\alpha_{n}^{-1}(\supp(\eta))}\|\pi(g)\eta(g^{-1} n)-\eta(n)\|\\
&=\varepsilon(1+K),
\end{split}\]
where $K=\max_{g\in\alpha_{n}^{-1}(\supp(\eta))}\|\pi(g)\eta(g^{-1} n)-\eta(n)\|$ is constant.
\end{proof}
Define a representation $\pi'$ of $H$ on the field by
\[\pi'(k)[h,n]:=[h,n\cdot k],\]
where $J_{H}(n)=t(k)$ (note that $J_{H}(n)=J_{H}(g\cdot n)$ for all $g\in G_{J(n)}$). Obviously, $\pi'$ is unitary.
Continuity of $\pi'$ easily follows from continuity of the right $H$ action on $N$. We formally denote this representation by $\HH\otimes_{G}N$. Define the map $\Theta_{N}:\RRep(G)\to\RRep(H)$ by
\[ [\pi,\HH,\Delta]\mapsto[\pi',\HH\otimes_{G}N,(J_{N}^{*}\Delta_{\HH})_{G}].\]
Let $[N^{-1}]:H\to G$ be the inverse of $[N]:G\to H$, i.e. $N\otimes_{K} N^{-1}\simeq G$, where $U(G)$ is $G$ itself seen as the unit $G-G$ bundle, and $N^{-1}\otimes_{G} N\simeq U(K)$. One easily shows that
\[\Theta_{N^{-1}}:[\pi,\HH,\Delta]\mapsto[\pi',\HH\otimes_{G}N^{-1},(J_{N^{-1}}^{*}\Delta_{\HH})_{G}]\]
is the inverse of $\Theta_{N}$.

It is obvious how to extend these maps to arrows.
\end{proof}

\subsection{Representation rings  and K-theory of a groupoid}
Suppose $G\rra M$ is a continuous groupoid and $M/G$ is compact.
\begin{definition}
The set of isomorphism classes of finite-dimensional continuous unitary representations of $G$, endowed with $\oplus$ and
$\otimes$ form a unital semi-ring. Applying the Grothendieck construction one obtains the \textbf{representation ring of} $G$,
denoted by $\RC_{f}(G)$. 
Denote the subring of locally trivial representations (projective Hilbert $C_{0}(M)$-modules) by $\RC(G)$. 
\end{definition}
\begin{example}
Suppose $M$ is a compact space. Consider the groupoid $M\rra M$. 
By definition one has $K_{0}(M)=\RC(M)$.
\end{example}
\begin{example}
Suppose $H$ is a compact group. Consider the groupoid $G:=H\rra pt$. Then $\RC_{f}(G)=\RC(G)$ equals the usual represenation ring $\RC(H)$ of $H$.
\end{example}
\begin{example}
Suppose $M$ is a locally compact space. Then for the pair groupoid $M\times M\rra M$ one sees that $\RC(M\times M)=\RC_{f}(M\times M)\simeq\ZZ$
generated by the trivial representation.
\end{example}
\begin{example}
One easily sees that Morita equivalent groupoids have isomorphic representation rings (as a corrolary of Theorem \ref{morrep}).
Hence, for a group $H$ and a principal $H$-bundle $P\to M$ one has 
\[\RC_{f}(P\times_{H}P)\simeq\RC(P\times_{H}P)\simeq \RC(H)\simeq \RC_{f}(H),\]
which generalizes the previous example.
\end{example}

Suppose $s,t:G\to M$ are open maps. Recall that the orbit relation of a groupoid $G\rra M$ is denoted by $R_{G}:=t\times s(G)$.
\begin{lemma}\label{modstruc}
The representation ring $\RC_{f}(G)$ is a $\RC_{f}(R_{G})$-module via the inclusion $\RC_{f}(R_{G})\to\RC_{f}(G)$ given by
\[\pi_{G}(g):=\pi_{R_{G}}(t(g),s(g)).\]
Analogously, $\RC(G)$ is a $\RC(R_{G})$-module.
\end{lemma}
\begin{example}
Suppose $s:G\to M$ is a continuous family of groups. Then $\RC_{f}(G)$ is a $\RC_{f}(M)$-module and
$\RC(G)$ is a $K_{0}(M)$-module (cf.\ Remark \ref{twist}). 
\end{example}

For proper groupoids the representation ring relates as follows to the K-theory of the reduced $C^{*}$-algebra of the groupoid. This was proved in \cite{tuxul}
in more general setting. We give a summary of their proof.
\begin{theorem}(\cite{tuxul})\label{txlthm}
If $G\rra M$ is proper, $M/G$ compact and $C^{*}_{r}(G)\otimes\mathcal{\mathbb{H}}$ has an approximate unit consisting of projections, then
$K_{0}(C^{*}_{r}(G))\simeq \RC(G)$.
\end{theorem}
\begin{proof}(sketch) One first shows, that $C^{*}_{r}(G)$ equals the $C^{*}$-algebra of compact
operators on $\hat{L}^{2}_{s}(G)\otimes\mathbb{H}$ made $G$-equivariant by an averaging construction. 
If $C^{*}_{r}(G)\otimes\mathcal{K}(\mathbb{H})$ has an approximate unit consisting of projections, then
$K_{0}(C^{*}_{r}(G))$ is obtained from the semi-ring generated by projections in $C^{*}_{r}(G)\otimes\mathbb{H}$, i.e. projections of
$\hat{L}^{2}_{s}(G)\otimes\mathbb{H}$. But these correspond precisely to locally trivial unitary representations of $G$ according to Proposition
\ref{tuxulprop} and the Serre-Swan theorem.
\end{proof}
\begin{example}
Suppose $M$ is a compact space. Then $C^{*}(M\rra M)=C_{0}(M)$, and $K_{0}(C_{0}(M))=K^{0}(M)=\RC(M\rra M)$.
Also, for the pair groupoid one can show $C^{*}_{r}(M\times M\rra M)\simeq\mathcal{K}(L^{2}(M))$ (cf.\ \cite{landsman2}) and
hence $K_{0}(C^{*}_{r}(M\times M\rra M))\simeq K_{0}(\mathcal{K}(L^{2}(M)))\simeq\ZZ\simeq \RC(M\times M\rra M)$.
\end{example}

\begin{example}
For a compact group $H$ it is well-known that $C^{*}_{r}(H)\simeq \hat{\bigoplus}_{(\pi,V)\in\hat{H}}M_{\dim(V)}(\CC)$ (with closure in the right norm, cf.\ \cite{landsman2}) and hence
\begin{align*}
K_{0}(C^{*}_{r}(H))&\simeq K_{0}(\hat{\bigoplus}_{(\pi,V)\in\hat{H}}M_{\dim(V)}(\CC))\\
&\simeq \bigoplus_{(\pi,V)\in\hat{H}}K_{0}(M_{\dim(V)}(\CC))\\
&\simeq \bigoplus_{(\pi,V)\in\hat{H}}\ZZ\simeq \RC(H).
\end{align*}
Theorem \ref{txlthm} generalizes this statement to proper groupoids (satisfying the mentioned condition).

For a principal $H$-bundle $P\to M$ one can prove $C^{*}_{r}(P\times_{H}P\rra M)\simeq C^{*}_{r}(H)\otimes\mathcal{K}(L^{2}(M))$, hence
\begin{align*}
K_{0}(C^{*}_{r}(P\times_{H}P\rra M))&\simeq K_{0}(C^{*}_{r}(H)\otimes\mathcal{K}(L^{2}(M)))\\
&\simeq K_{0}(C^{*}_{r}(H))\\
&\simeq \RC(H)\simeq\RC(P\times_{H}P\rra M),
\end{align*}
by stablity of $K$-theory.
\end{example}

\section{The groupoid convolution Banach *-category}\label{grconv}
In \cite{ren} Reneault established a bijective correspondence between representations of groupoids on measurable fields of Hilbert spaces and the 
non-degenerate bounded representations of the Banach *-algebra $L^{1}(G)$, generalizing the analogous statement for groups. In this section we shall 
prove a different generalization suitable for continuous 
representations of groupoids. We give a bijective correspondence between continuous representations of groupoids on continuous fields of Hilbert 
spaces and continuous representations on continuous fields of Hibert spaces of the Banach *-category $\hat{L}^{1}(G)$.

\subsection{Fell bundles over groupoids and continuous $C^{*}$-categories.}
First we need some terminology. We discuss the relation between continuous Fell bundles over groupoids (cf.\ \cite{yam, muh, kum}) and
Banach $*$-categories and $C^{*}$-categories (cf.\ \cite{doro}).

A \textbf{(lower semi-)continuous Fell bundle over a groupoid} $G$ is a (lower \mbox{semi-)}continuous field of Banach spaces $(\{\BB_{g}\}_{g\in G},\Delta)$ over
$G$ endowed with
an associative bilinear product 
\[\BB_{g}\times\BB_{h}\to\BB_{g h}, (P,Q)\mapsto P Q\]
whenever $(g,h)\in G^{(2)}$ and an anti-linear involution
\[\BB_{g}\to\BB_{g^{-1}}, P\mapsto P^{*}\]
satisfying the following conditions for all $(g,h)\in G^{(2)}$ and $(P,Q)\in \BB_{g}\times\BB_{h}$
\begin{itemize}
\item[(i)] $\|P Q\|\leq\|P\|\|Q\|$;
\item[(ii)] $\|P^{*}P\|=\|P\|^{2}$;
\item[(iii)] $(P Q)^{*}=Q^{*}P^{*}$;
\item[(iv)] $P^{*}P$ is a positive element of $\BB_{1_{s(g)}}$;
\item[(v)] the image of the multiplication $\BB_{g}\times\BB_{h}\to\BB_{g h}, (P,Q)\mapsto P Q$ is dense;
\item[(vi)] multiplication $m^{*}\BB\to\BB$ and involution $\BB\to\BB$ are continuous maps of fields of Banach spaces.
\end{itemize}
where $\BB$ denotes the total space of $(\BB_{m\in M},\Delta)$ endowed with the topology given by $\Delta$ and $m^{*}\BB$ the pullback of the field $\BB$ over $G$ along $m: G^{(2)}\to G$.

\begin{example} Our main example will be the following. Suppose $G\rra M$ is a continuous groupoid with open $s,t:G\to M$ and
\mbox{$R_{G}:=(t\times s)(G)$}. Let $(\{\HH_{m}\}_{m\in M},\Delta_{\HH})$ be a continuous field of Hilbert spaces over $M$.
Consider the lower semi-continuous field of Banach spaces over $R_{G}$ whose fiber at $(n,m)$ is given by the bounded linear operators
$\HH_{m}\to\HH_{n}$, i.e. $\BB_{(n,m)}:=\BB(\HH_{n},\HH_{m})$. This field was already introduced in Section \ref{opno}.
\begin{lemma}\label{BBlem}
The lower semi-continuous field of Banach spaces $(\{\BB_{(n,m)}\}_{(n,m)\in R_{G}},\Delta_{\BB})$ is a lower semi-continuous Fell bundle over $R_{G}$.
\end{lemma}
\begin{proof}
The continuity of the composition was proven in the proof of Lemma \ref{oppietoppie2}.
Note that $\pi(g):\HH_{s(g)}\to \HH_{t(g)}$ is an isomorphism of Hilbert spaces. Hence, the properties (i), (ii),
(iii),(iv) and (v) follow from the fact that these are true for $\BB(\HH)$, where
$\HH\simeq\HH_{t(g)}\simeq\HH_{s(g)}$.
\end{proof}
We sometimes write $\BB(\HH,\HH)$ for this lower semi-continuous Fell bundle over $R_{G}$.
\end{example}

A (lower semi-)continuous Fell bundle $A$ over a continuous equivalence relation $R\subset M\times M$ on $M$ is a \textbf{(lower semi-)continuous $C^{*}$-category over} $M$. Leaving out the $C^{*}$-norm equality (ii) we speak of a \textbf{(lower semi-)continuous Banach $*$-category}. Note that it is indeed a category with well-defined source and target maps $s,t:A\to M$.

\begin{example}
Let $G\rra M$ be a locally compact groupoid endowed with a Haar system $\{\lambda^{m}\}_{m\in M}$. 
Suppose there exist a continuous families measures $\{\lambda_{m}^{n}\}_{(n,m)\in R_{G}}$ on $G$ and $\{\mu_{m}\}_{m\in M}$ on $M$ such that
\[\lambda^{n}=\int_{m\in s(G^{n})}\lambda_{m}^{n}\mu_{n}(dm).\]

Consider the continuous field of Banach spaces $(\hat{L}^{1}(G),\Delta^{1}(G)):=(\hat{L}^{1}_{t\times s}(G),\Delta^{1}_{t\times s}(G))$, cf.\
Example \ref{exfield}.
\begin{lemma}
$(L^{1}(G),\Delta^{1}(G))$ is a continuous Banach $*$-category over $M$, where the multiplication map $\hat{L}^{1}(G)^{(2)}\to \hat{L}^{1}(G)$ is the continuous extension of
\[f * f'(g):=\int_{h\in G^{m}_{k}}f(g h^{-1})f'(h)\lambda_{k}^{m}(d h),\]
for all $f\in C_{c}(G_{m}^{n})$ and $f'\in C_{c}(G_{k}^{m})$.
\end{lemma}
\begin{proof}(sketch) One, indeed easily checks that $\|f*g\|(n,k)<\|f\|(n,m)\|g\|(m,k)$, so this extension is well-defined.
\end{proof}
\end{example}

\begin{definition}
A \textbf{strongly continuous representation $(\HH,\Delta,L)$ of a continuous Banach $*$-category} $A$ over a space $M$ on a continuous field of Hilbert spaces $(\HH,\Delta)$ over $M$ is a
continuous $*$-homomorphism 
\[L:A\to\BB(\HH,\HH),\]
such that $a\mapsto L(a)\xi(s(a))$ is continuous $A\to\HH$ for every $\xi\in\Delta$.
\end{definition}
One has analogous definitions for weakly continuous representations and representations continuous in the operator norm.
\begin{definition}
A representation $(\HH,\Delta, L)$ of a Banach $*$-category $A$ is non-degenerate if $\overline{L(A)\HH}=\HH$.
\end{definition}

\subsection{Representations of $G$ versus representations of $\hat{L}^{1}(G)$.}
In this section we again need Dirac sequences but in a different way.
Suppose $G\rra M$ allows Dirac sequences $\{(\delta_{k}^{g})_{k\in\NN}\}_{g\in G}$ for the Haar system $\{\lambda_{m}^{n}\}_{(n,m)\in R_{G}}$ in the sense that
\begin{itemize}
\item[i)] $\delta_{k}^{g}\geq 0$ on $G_{s(g)}^{t(g)}$,
\item[ii)] $\int_{g'\in G_{m}^{n}}\delta_{k}^{g}(g')\lambda_{s(g)}^{t(g)}(d g')=1$ for all $k\in\NN$,
\item[iii)] For every open neighborhood $U\subset G^{n}_{m}$ of $g$ and every $\varepsilon>0$ there is an $N\in\NN$ such that for $k>N$
\[\int_{g\in U^{c}}\delta^{g}_{k}(g')\lambda^{n}_{m}(d g')<\varepsilon.\]
\end{itemize}

\begin{lemma}\label{heen}
If $(\HH,\Delta,\pi)$ is a continuous unitary representation of $G\rra M$, then $L_{\pi}:\hat{L}^{1}(G)\to\BB(\HH,\HH)$ given by
\[f\mapsto\left((n,m)\mapsto\int_{G_{m}^{n}}f(g)\pi(g)\lambda_{m}^{n}(dg)\right)\]
is a non-degenerate strongly continuous representation of $(\hat{L}^{1}(G),\Delta{^1}(G))$ as a continuous Banach $*$-category on the continuous field of Hilbert spaces $(\HH,\Delta)$.
\end{lemma}
\begin{proof}
By the properties of the Bochner integral one has
\[\begin{array}{ll}
\|\pi(f)\|(n,m)&=\|\int_{G_{m}^{n}}f(g)\pi(g)\lambda^{n}_{m}\|\\
&\leq\int_{G_{m}^{n}}|f(g)| \|\pi(g)\|\lambda_{m}^{n}(d g)\\
&=\int_{G_{m}^{n}}|f(g)|\lambda_{m}^{n}(d g)
\end{array}\]
(Note that therefore, $\|\pi(f)\|=\sup_{(n,m)\in R}\|\pi(f)\|(n,m)\leq\|f\|_{\hat{L}^{1}(G)}).$

We now prove that $L_{\pi}$ is a $*$-homomorphism. Suppose $f\in C_{c}(G_{m}^{k})$ and $f'\in C_{c}(G_{k}^{n})$, then
\begin{align*}
L_{\pi}(f * f')&=\int_{g\in G_{m}^{n}}(f * f')(g)\pi(g)\lambda_{m}^{n}(d g)\\
&=\int_{g\in G_{m}^{n}}\int_{h\in G_{m}^{k}}f(g h^{-1})f'(h)\lambda_{m}^{k}(d h)\pi(g)\lambda_{m}^{n}(d g)\\
&=\int_{g\in G_{k}^{n}}f(g)\pi(g)\lambda_{k}^{n}(d g)\int_{h\in G_{m}^{k}}f'(h)\pi(h)\lambda_{k}^{m}(d h)\\
&=L_{\pi}(f)L_{\pi}(f'),
\end{align*}
by invariance of the Haar system.

Suppose $f\in C_{c}(G_{m}^{n})$ is given. Suppose $F\in C_{c}(G)$ satisfies $F|_{G_{m}^{n}}=f$.
Note that 
\[\sup_{(n',m')\in R_{G}}\int_{G_{m}^{n}}\|F(g)\pi(g)\xi(s(g))\|\lambda_{m}^{n}(d g)\leq\|F\|_{\hat{L}^{1}(G)}\max_{m\in s(\supp(F))}\|\xi(m)\|.\]
Using this, one easily proves that $L_{\pi}$ is strongly continuous.

The representation $\pi_{L}$ is non-degenerate, since for any $m\in M$ and $h\in\HH_{m}$
\[\lim_{k\to\infty}\|h-L(\delta^{1_{m}}_{k})h\|=\lim_{k\to\infty}\|h-\int_{g\in G_{m}^{n}}\delta^{1_{m}}_{k}(g)\pi(g)h\lambda_{m}^{n}(d g)\|=0.\]
This finishes the proof.
\end{proof}

For $f\in C_{c}(G)$, $m\in M$ and $g,g'\in G^{m}$, we shall use the notation $f^{g}(g'):=(\pi_{L}(g)f)(g')=f(g^{-1}g')$.
\begin{lemma}\label{terug}
If $(\HH,\Delta, L)$ is a strongly continuous non-degenerate representation of \linebreak $(\hat{L}^{1}(G),\Delta^{1}(G))$, then
\[\pi(g)(L(f)h):=L(f^{g})h\]
defines a continuous unitary representation of $G$ on the continuous field of Hilbert spaces $(\HH,\Delta)$.
\end{lemma}
\begin{proof}
By non-degeneracy of $L$, the above formula defines $\pi_{L}$ on a dense set. 
It extends to the whole of $\HH$, since for all $g\in G$ and $h\in\HH_{s(g)}$ one has 
\[\begin{array}{ll}
\|\pi(g)h\|&=\lim_{k\to\infty}\|L(\delta_{g}^{k})h\|\\
&\leq \lim_{k\to\infty}B\|\delta_{g}^{k}\|\|h\|\\
&=B\|h\|,
\end{array}\]
for a constant $B\in\RR\geq 0$.

This is well-defined. Indeed, suppose $L(f)h=L(f')h'$ for
$f\in L^{1}(G_{m}^{n})$, $f'\in L^{1}(G_{m'}^{n})$, $h\in\HH_{m}$ and $h'\in\HH_{m'}$.
Let $\delta_{k}^{g}$ denote the translation of $\delta_{k}^{m}$ along $g$. One easily checks that
\[\|\delta^{g}_{k}*f-f^{g}\|\to 0\]
when $k\to\infty$.
Then one has for all $k\in\NN$:
\[\begin{array}{ll}
\|L((f'^){g})h'-L(f^{g})h\|&\leq \|L((f')^{g})h'-L(\delta_{k}^{g}*f')h'\|\\
&\|L(\delta_{k}^{g}*f')h'-L(\delta_{k}^{g}*f)h\|+\|L(\delta_{k}^{g}*f)h'-L(f^{g})h\|\\
&\leq B \|(f')^{g}-\delta^{g}_{k}*f'\|\|h'\|+\|L(\delta_{k}^{g})(L(f')h'-L(f)h)\|\\
&+B \|(f)^{g}-\delta^{g}_{k}*f\|\|h\|
\end{array}\]
The second term is zero and the first and the last term go to zero as $k\to\infty$, hence
$L((f'^){g})h'=L(f^{g})h$.

$\pi$ is a homomorphism. Indeed, for $(g,g')\in G^{(2)}$, $f\in L^{1}(G_{m}^{s(k)})$ and $h\in\HH_{m}$ one has
\[\begin{array}{ll}
\pi(g g')(L(f)h)&=L(f^{g g'})h\\
&=L((f^{g'})^{g})h\\
&=\pi(g)L(f^{g'})h\\
&=\pi(g)\pi(g')(L(f)h).
\end{array}\]

Furthermore, the following computation shows that $\pi(g)^{*}=\pi(g^{-1})$:
\[\begin{array}{ll}
\la\pi(g)^{*}L(f)h,L(f')h'\ra&=\la h,L(f)^{*}\pi(g)L(f')h'\ra\\
&=\la h,L(f^{*})L((f')^{g})h')h'\ra\\
&=\la h,L(f^{*}*(f')^{g})h'\ra\\
&=\la h,L((f^{g^{-1}})^{*}*f')h'\ra\\
&=\la h,L(f^{g^{-1}})^{*}L(f')h'\ra\\
&=\la L(f^{g^{-1}})h,L(f')h'\ra\\
&=\la \pi(g^{-1})L(f)h,L(f')h'\ra,
\end{array}\]
where the fourth step follows from equivariance of the Haar system and the fact that
\[(f^{g})^{*}(g')=f^{*}(g' g^{-1}).\]

The continuity of $\pi$ follows from the fact that for any $F\in C_{c}(G)$, representing a section of $\hat{L}^{1}(G)\to R_{G}$,
and any $\xi\in\Delta$, the section $m\mapsto L(F)(m,m)\xi(m)$ is again in $\Delta$ and that $g\mapsto F(s(g),s(g))^{g}$ is continuous,
cf.\ Lemma \ref{regrep}.
\end{proof}

\begin{theorem}
The correspondence $\pi\mapsto L_{\pi}$ is a bijection between the set of continuous unitary representations of $G$ and the set of strongly continuous
representations of $(\hat{L}^{1}(G),\Delta^{1}(G))$. 
\end{theorem}
\begin{proof}
The inverse correspondence is given by Lemma \ref{terug}, which we denote by $L\mapsto \pi^{L}$ (not to be confused with the left regular representation $\pi_{L}$). Given a continuous unitary representation $\pi$ of $G$, we compute
\[\begin{array}{ll}
\pi^{(L^{\pi})}(g)(L^{\pi}(f)h)&=L^{\pi}(f^{g})h\\
&=\int_{g'\in G_{m}^{n}}f(g^{-1}g')\pi(g')h\lambda_{m}^{n}(d g')\\
&=\int_{g'\in G_{m}^{n}}f(g')\pi(g)\pi(g')h\lambda_{m}^{p}(d g')\\
&=\pi(g)(L^{\pi}(f)h).
\end{array}\]
Conversely, suppose a non-degenerate strongly continuous representation $L$ of $\hat{L}^{1}(G)$ is given. Then we have
\[\begin{array}{ll}
L^{(\pi^{L})}L(f')h&=\int_{g\in G_{m}^{n}}f(g)\pi^{L}(g)L(f')h\lambda_{m}^{n}(d g)\\
&=\int_{g\in G_{m}^{n}}f(g)L((f')^{g})h\lambda_{m}^{n}(d g)\\
&=L(\int_{g\in G_{m}^{n}}f(g)(f')^{g}\lambda_{m}^{n}(d g))h\\
&=L(f*f')h=L(f)(L(f')h),
\end{array}\]
which finishes the proof.
\end{proof}

\bibliographystyle{plain}

\end{document}